\documentclass[a4paper,12pt]{article}

\usepackage{latexsym}
\usepackage{amssymb}
\usepackage{theorem}
\usepackage{amsmath}
\usepackage{amscd}
\usepackage[dvips]{graphicx}

\newtheorem{theorem}{Theorem}[section]
\newtheorem{corollary}[theorem]{Corollary}
\newtheorem{lemma}[theorem]{Lemma}
\newtheorem{example}[theorem]{Example}
\newtheorem{proposition}[theorem]{Proposition}
\newtheorem{remark}[theorem]{Remark}
\newtheorem{definition}[theorem]{Definition}

\newcommand{\demo}{\par\noindent{\it Proof. \/}\ }
\newcommand{\enD}{\hfill $\Box$\vspace{3truemm} \par}
\newcommand{\R}{\mathbb{R}}

\newcommand{\bn}{\mbox{\boldmath $n$}}
\newcommand{\bt}{\mbox{\boldmath $t$}}

\newcommand{\bb}{\mbox{\boldmath $b$}}
\newcommand{\ba}{\mbox{\boldmath $a$}}
\newcommand{\bc}{\mbox{\boldmath $c$}}
\newcommand{\bd}{\mbox{\boldmath $d$}}

\newcommand{\be}{\mbox{\boldmath $e$}}
\newcommand{\bmu}{\mbox{\boldmath $\mu$}}

\begin{document}

\title{Bertrand and Mannheim curves of framed curves in the $4$-dimensional Euclidean space}

\author{
Shun'ichi Honda, Masatomo Takahashi and Haiou Yu}

\date{\today}

\maketitle

\begin{abstract}
A Bertrand curve in the $4$-dimensional Euclidean space is a space curve whose first normal line is the same as the first normal line of another curve.
On the other hand, a Mannheim curve in the $4$-dimensional Euclidean space is a space curve whose first normal line is the same as the second or third normal line of another curve.
By definitions, another curve is a parallel curve with respect to the direction of the first normal vector.
As smooth curves with singular points, we consider framed curves in the Euclidean space.
Then we define and investigate Bertrand and Mannheim curves of framed curves.
We give necessary and sufficient conditions of Bertrand and Mannheim curves of both regular and framed curves.
It is well-known that the Bertrand curves of regular curves do not exist under a condition.
However, even if regular curves, Bertrand curves exist as framed curves.
\end{abstract}

\renewcommand{\thefootnote}{\fnsymbol{footnote}}

\section{Introduction}

Bertrand and Mannheim curves are classical objects in differential geometry (\cite{Aminov,Banchoff-Lovett,Berger-Gostiaux,Bertrand,doCarmo,Izumiya-Takeuchi1,Kuhnel,Liu-Wang,Matsuda-Yorozu1,Matsuda-Yorozu2,Pears,Struik}).
A Bertrand curve in $4$-dimensional Euclidean space is a space curve whose first normal line is the same as the first normal line of another curve (cf. \cite{Matsuda-Yorozu1,Pears}).
On the other hand, a Mannheim curve in $4$-dimensional Euclidean space is a space curve whose first normal line is the same as the second or third normal line of another curve.
By definitions, another curve is a parallel curve with respect to the direction of the first normal vector.
In order to define the Frenet frame, a condition (a non-degenerate condition) is needed.
In general, the parallel curve does not satisfy these conditions.

It is well-known that the Bertrand curves of regular curves do not exist under a condition in \cite{Matsuda-Yorozu1,Pears}.
In \cite{Matsuda-Yorozu2}, they consider the condition of the Mannheim curves of regular curves.
We clarify the existence conditions of Mannheim curves of regular curves in \S 2.

As smooth curves with singular points, we introduced the notion of framed curves in the Euclidean space in \cite{Honda-Takahashi1}.
Then we define Bertrand and Mannheim curves of framed curves in \S 3 and \S 4, respectively.
We give necessary and sufficient conditions of the Bertrand and Mannheim curves of framed curves, respectively.
Even if regular curves, Bertrand curves exist as framed curves in \S 5.
The basic results on the singularity theory see \cite{Arnold1,Arnold2,Bruce-Giblin,Ishikawa-book,Izumiya-book}.

All maps and manifolds considered in this paper are differentiable of class $C^{\infty}$.

\bigskip
\noindent
{\bf Acknowledgement}.
The second author was partially supported by JSPS KAKENHI Grant Number JP 20K03573.
The third author was partially supported by Natural Science Foundation of Jilin Province of China Grant Number JC 20210101153.

\section{Preliminaries}

Let $\R^4$ be the $4$-dimensional Euclidean space equipped with the inner product $\ba \cdot \bb = a_1 b_1 + a_2 b_2 + a_3 b_3+a_4 b_4$,
where $\ba = (a_1, a_2, a_3, a_4)$ and $\bb = (b_1, b_2, b_3, b_4) \in \R^4$.
The norm of $\ba$ is given by $\vert \ba \vert = \sqrt{\ba \cdot \ba}$ and the vector product is given by
$$
\ba \times \bb \times \bc={\rm det}
\left(
\begin{array}{cccc}
a_1 & a_2 & a_3 & a_4\\
b_1 & b_2 & b_3 & b_4\\
c_1 & c_2 & c_3 & c_4\\
\be_1 & \be_2 & \be_3 & \be_4
\end{array}
\right),
$$
where $\{\be_1, \be_2, \be_3, \be_4\}$ is the canonical basis of $\R^4$ and $\bc=(c_1,c_2,c_3,c_4) \in \R^4$.
Let $S^3$ be the unit sphere in $\R^4$, that is, $S^3=\{\ba \in \R^4 \ |\ |\ba|=1\}$.
We denote the $6$-dimensional smooth manifold $\{(\ba,\bb,\bc) \in S^3 \times S^3 \times S^3 \ |\ \ba \cdot \bb=\ba \cdot \bc=\bb \cdot \bc=0\}$ by $\Delta$.

Note that for $(\ba,\bb,\bc) \in \Delta$, we denote $\ba \times \bb \times \bc=\bd$, then
$$
\bd \times \ba \times \bb=-\bc, \ \bc \times \bd \times \ba=\bb, \ \bb \times \bc \times \bd=-\ba.
$$

We quickly review the theories of Bertrand curves, Mannheim curves of regular curves and framed curves.

\subsection{Regular space curves in $\R^4$}

Let $I$ be an interval of $\R$ and let $\gamma:I \to \R^4$ be a regular space curve, that is, $\dot{\gamma}(t) \not=0$ for all $t \in I$, where $\dot{\gamma}(t)=(d\gamma/dt)(t)$.
We say that $\gamma$ is {\it non-degenerate}, or $\gamma$ satisfies the {\it non-degenerate condition} if $\dot{\gamma}(t) \times \ddot{\gamma}(t) \times \dddot{\gamma}(t) \not=0$ for all $t \in I$.
We take the arc-length parameter $s$, that is, $|\gamma'(s)|=1$ for all $s$, where $\gamma'(s)=(d\gamma/ds)(s)$.
Then the tangent, the first, second and third normal maps $\bt,\bn_1,\bn_2,\bn_3:I \to S^3$ are given by
\begin{eqnarray*}
&&\bt(s)=\gamma'(s), \ \bn_1(s)=\frac{\gamma''(s)}{|\gamma''(s)|}=\frac{\bt'(s)}{|\bt'(s)|}, \\ &&\bn_2(s)=\frac{\bn'_1(s)+|\bt'(s)|\bt(s)}{|\bn'_1(s)+|\bt'(s)|\bt(s)|}, \ \bn_3(s)=\bt(s) \times \bn_1(s) \times \bn_2(s)
\end{eqnarray*}
under the non-degenerate condition $\gamma'(s) \times \gamma''(s) \times \gamma'''(s) \not=0$ for all $s$.
Then $\{\bt(s),\bn_1(s),\bn_2(s),\bn_3(s)\}$ is a moving frame of $\gamma(s)$ and we have the Frenet-Serret formula (cf. \cite{Aminov,Lee}):
$$
\left(
\begin{array}{c}
\bt'(s)\\
\bn'_1(s)\\
\bn'_2(s)\\
\bn'_3(s)
\end{array}
\right)
=
\left(
\begin{array}{cccc}
0&\kappa_1(s)&0&0\\
-\kappa_1(s)&0&\kappa_2(s)&0\\
0&-\kappa_2(s)&0&\kappa_3(s)\\
0&0&-\kappa_3(s)&0
\end{array}
\right)
\left(
\begin{array}{c}
\bt(s)\\
\bn_1(s)\\
\bn_2(s)\\
\bn_3(s)
\end{array}
\right),
$$
where
\begin{eqnarray*}
\kappa_1(s)=|\gamma''(s)|, \ \kappa_2(s)=\frac{|\gamma'(s) \times \gamma''(s) \times \gamma'''(s)|}{\kappa_1^2(s)}, \ \kappa_3(s)=\frac{{\rm det}(\gamma'(s),\gamma''(s),\gamma'''(s),\gamma^{(4)}(s))}{\kappa^3_1(s)\kappa^2_2(s)}
\end{eqnarray*}
are the first, second and third curvatures of $\gamma(s)$, respectively.
If we take a general parameter $t$, then the tangent, the first, second and third normal maps $\bt, \bn_1,\bn_2,\bn_3:I \to S^3$ are given by
\begin{eqnarray*}
&&\bt(t)=\frac{\dot{\gamma}(t)}{|\dot{\gamma}(t)|}, \ \bn_1(t)=\frac{\dot{\bt}(t)}{|\dot{\bt}(t)|}=\frac{|\dot{\gamma}(t)|^2\ddot{\gamma}(t)-(\dot{\gamma}(t) \cdot \ddot{\gamma}(t))\dot{\gamma}(t)}{||\dot{\gamma}(t)|^2\ddot{\gamma}(t)-(\dot{\gamma}(t) \cdot \ddot{\gamma}(t))\dot{\gamma}(t)|}, \\ &&\bn_2(t)=-\bn_3(t) \times \bt(t) \times \bn_1(t), \ \bn_3(t)=\frac{\dot{\gamma}(t) \times \ddot{\gamma}(t) \times \dddot{\gamma}(t)}{|\dot{\gamma}(t) \times \ddot{\gamma}(t) \times \dddot{\gamma}(t)|}
\end{eqnarray*}
under the non-degenerate condition $\dot{\gamma}(t) \times \ddot{\gamma}(t) \times \dddot{\gamma}(t) \not=0$ for all $t \in I$.
Then $\{\bt(t),\bn_1(t),\bn_2(t),\bn_3(t)\}$ is a moving frame of $\gamma(t)$ and we have the Frenet-Serret formula:
$$
\left(
\begin{array}{c}
\dot{\bt}(t)\\
\dot{\bn}_1(t)\\
\dot{\bn}_2(t)\\
\dot{\bn}_3(t)
\end{array}
\right)
=
\left(
\begin{array}{cccc}
0&|\dot{\gamma}(t)|\kappa_1(t)&0&0\\
-|\dot{\gamma}(t)|\kappa_1(t)&0&|\dot{\gamma}(t)|\kappa_2(t)&0\\
0&-|\dot{\gamma}(t)|\kappa_2(t)&0&|\dot{\gamma}(t)|\kappa_3(t)\\
0&0&-|\dot{\gamma}(t)|\kappa_3(t)&0
\end{array}
\right)
\left(
\begin{array}{c}
\bt(t)\\
\bn_1(t)\\
\bn_2(t)\\
\bn_3(t)
\end{array}
\right),
$$
where
\begin{eqnarray*}
\kappa_1(t) &=& \frac{\sqrt{|\dot{\gamma}(t)|^2|\ddot{\gamma}(t)|^2-(\dot{\gamma}(t) \cdot \ddot{\gamma}(t))^2}}{|\dot{\gamma}(t)|^3}, \\
\kappa_2(t) &=& \frac{|\dot{\gamma}(t) \times \ddot{\gamma}(t) \times \dddot{\gamma}(t)|}{|\dot{\gamma}(t)|^6 \kappa^2_1(t)}, \\
\kappa_3(t) &=& \frac{{\rm det}(\dot{\gamma}(t),\ddot{\gamma}(t),\dddot{\gamma}(t),\gamma^{(4)}(t))}{|\dot{\gamma}(t)|^{10} \kappa^3_1(t) \kappa_2^2(t)}.
\end{eqnarray*}
\par

Note that in order to define the moving frame $\{\bt(t),\bn_1(t),\bn_2(t),\bn_3(t)\}$ of $\gamma(t)$ and the curvatures, we assume that $\gamma$ is not only regular, but also non-degenerate.
\par
It is known that $|\dot{\gamma}(t)|^2|\ddot{\gamma}(t)|^2-(\dot{\gamma}(t) \cdot \ddot{\gamma}(t))^2 \geq 0$ for all $t\in I$.
If $|\dot{\gamma}(t_0)|^2|\ddot{\gamma}(t_0)|^2-(\dot{\gamma}(t_0) \cdot \ddot{\gamma}(t_0))^2=0$ for $t_0 \in I$,
we have $\cos\angle(\dot{\gamma}(t_0),\ddot{\gamma}(t_0))=\pm1$.
By $\gamma$ is regular, there exists a constant $a\in\R$ such that $\ddot{\gamma}(t_0)=a\dot{\gamma}(t_0)$.
It follows that $\dot{\gamma}(t_0) \times \ddot{\gamma}(t_0) \times \dddot{\gamma}(t_0)=0$.
Thus, if $\dot{\gamma}(t) \times \ddot{\gamma}(t) \times \dddot{\gamma}(t)\neq0$ for all $t\in I$,
then $|\dot{\gamma}(t)|^2|\ddot{\gamma}(t)|^2-(\dot{\gamma}(t) \cdot \ddot{\gamma}(t))^2>0$ for all $t\in I$.
Namely, if $\gamma$ is non-degenerate, then the curvatures $\kappa_1$ and $\kappa_2$ are positive.
\subsection{Bertrand curves of regular space curves in $\R^4$}
Let $\gamma$ and $\overline{\gamma}:I \to \R^4$ be different non-degenerate curves.
We denote the moving frames of $\gamma$ and $\overline{\gamma}$ by $\{\bt,\bn_1,\bn_2,\bn_3\}$ and $\{\overline{\bt},\overline{\bn}_1,\overline{\bn}_2,\overline{\bn}_3\}$, respectively.

\begin{definition}\label{regular-Bertrand}{\rm
We say that $\gamma$ and $\overline{\gamma}$ are {\it Bertrand mates} ({\it first mates} or {\it $(\bn_1,\overline{\bn}_1)$-mates}) if there exists a smooth function $\lambda:I \to \R$ such that $\overline{\gamma}(t)=\gamma(t)+\lambda(t)\bn_1(t)$ and $\bn_1(t)=\pm \overline{\bn}_1(t)$ for all $t \in I$.
We also say that $\gamma:I \to \R^4$ is a {\it Bertrand curve} if there exists another non-degenerate curve $\overline{\gamma}:I \to \R^4$ such that $\gamma$ and $\overline{\gamma}$ are Bertrand mates.
}
\end{definition}

If $\gamma$ and $\overline{\gamma}$ are Bertrand mates, then the first normal line of $\gamma$ and the first normal line of $\overline{\gamma}$ are the same for each point.

By a parameter change, we may assume that $s$ is the arc-length parameter of $\gamma$.

\begin{lemma}\label{lambda-const-Bertrand}
Let $\gamma:I \to \R^4$ be non-degenerate with the arc-length parameter.
Under the notations in Definition \ref{regular-Bertrand}, if $\gamma$ and $\overline{\gamma}$ are Bertrand mates, then $\lambda$ is a non-zero constant.
\end{lemma}
\demo
By differentiating $\overline{\gamma}(s)=\gamma(s)+\lambda(s)\bn_1(s)$,
we have
$$
|\dot{\overline{\gamma}}(s)|\overline{\bt}(s)=(1-\lambda(s)\kappa_1(s))\bt(s)+\lambda'(s)\bn_1(s)+\lambda(s)\kappa_2(s)\bn_2(s).
$$
Since $\bn_1(s)= \pm \overline{\bn}_1(s)$, we have $\lambda'(s)=0$ for all $s \in I$. Therefore $\lambda$ is a constant.
If $\lambda=0$, then $\overline{\gamma}(s)=\gamma(s)$ for all $s \in I$.
Hence, $\lambda$ is a non-zero constant.
\enD
\begin{theorem}\label{regular-Bertrand-existence}
Let $\gamma:I \to \R^4$ be non-degenerate with the arc-length parameter.
Suppose that $\kappa_3(s_0) \not=0$ for a point $s_0 \in I$.
Then $\gamma$ is not a Bertrand curve.
\end{theorem}
\demo
Suppose that $\gamma$ is a Bertrand curve.
There exists another non-degenerate curve $\overline{\gamma}:I \to \R^4$ such that $\gamma$ and $\overline{\gamma}$ are Bertrand mates.
Note that the parameter $s$ is not the arc-length parameter of $\overline{\gamma}$.
By Lemma \ref{lambda-const-Bertrand}, $\overline{\gamma}(s)=\gamma(s)+\lambda \bn_1(s)$, where $\lambda$ is a non-zero constant and $\bn_1(s)=\pm \overline{\bn}_1(s)$.
By differentiating $\overline{\gamma}(s)=\gamma(s)+\lambda \bn_1(s)$, we have
$$
|\dot{\overline{\gamma}}(s)|\overline{\bt}(s)=(1-\lambda \kappa_1(s))\bt(s)+\lambda \kappa_2(s)\bn_2(s).
$$
There exists a smooth function $\theta:I \to \R$ such that $\overline{\bt}(s)=\cos \theta(s) \bt(s)+\sin \theta(s) \bn_2(s)$. That is,
$$
\cos \theta(s)=\frac{1-\lambda \kappa_1(s)}{|\dot{\overline{\gamma}}(s)|}, \ \sin \theta (s)=\frac{\lambda \kappa_2(s)}{|\dot{\overline{\gamma}}(s)|}.
$$
By differentiating $\overline{\bt}(s)=\cos \theta(s) \bt(s)+\sin \theta(s) \bn_2(s)$, we have
\begin{eqnarray*}
|\dot{\overline{\gamma}}(s)|\overline{\kappa_1}(s)\overline{\bn}_1(s)&=&-\theta'(s)\sin \theta(s) \bt(s)+(\kappa_1(s)\cos \theta(s)-\kappa_2(s)\sin \theta(s))\bn_1(s)\\
&&+\theta'(s)\cos \theta(s) \bn_2(s)+\kappa_3(s)\sin \theta(s)  \bn_3(s).
\end{eqnarray*}
Since $\bn_1(s)=\pm \overline{\bn}_1(s)$, ${\theta}(s)$ is a constant and $\kappa_3(s) \sin \theta =0$ for all $s \in I$.
However, it is the contradict to the fact that $\sin \theta \not=0$ and $\kappa_3(s_0) \not=0$.
Therefore, $\gamma$ is not a Bertrand curve.
\enD

\begin{remark}{\rm
If $\kappa_3(s)=0$ for all $s \in I$, then the trace of $\gamma$ is contained in a $3$-dimensional space.
Hence $\gamma$ is possible to be a Bertrand curve, see \cite{Honda-Takahashi2}.
}
\end{remark}

\subsection{Mannheim curves of regular space curves in $\R^4$}

Let $\gamma$ and $\overline{\gamma}:I \to \R^4$ be different non-degenerate curves.
We denote the moving frames of $\gamma$ and $\overline{\gamma}$ by $\{\bt,\bn_1,\bn_2,\bn_3\}$ and $\{\overline{\bt},\overline{\bn}_1,\overline{\bn}_2,\overline{\bn}_3\}$, respectively.

\begin{definition}\label{regular-Mannheim}{\rm
$(1)$ We say that $\gamma$ and $\overline{\gamma}$ are {\it second type of Mannheim mates} (briefly, {\it second mates} or {\it $(\bn_1,\overline{\bn}_2)$-mates}) if there exists a smooth function $\lambda:I \to \R$ such that $\overline{\gamma}(t)=\gamma(t)+\lambda(t)\bn_1(t)$ and $\bn_1(t)=\pm \overline{\bn}_2(t)$ for all $t \in I$.
\par
$(2)$ We say that $\gamma$ and $\overline{\gamma}$ are {\it third type of Mannheim mates} (briefly, {\it third mates} or {\it $(\bn_1,\overline{\bn}_3)$-mates}) if there exists a smooth function $\lambda:I \to \R$ such that $\overline{\gamma}(t)=\gamma(t)+\lambda(t)\bn_1(t)$ and $\bn_1(t)=\pm \overline{\bn}_3(t)$ for all $t \in I$.
\par
We also say that $\gamma:I \to \R^4$ is a {\it second type of Mannheim curve} (respectively, {\it third type of Mannheim curve}) if there exists another non-degenerate curve $\overline{\gamma}:I \to \R^4$ such that $\gamma$ and $\overline{\gamma}$ are second mates (respectively, third mates).
}
\end{definition}

If $\gamma$ and $\overline{\gamma}$ are second mates (respectively, third mates), then the first normal line of $\gamma$ and the second normal line (respectively, third normal line) of $\overline{\gamma}$ are the same for each point.

By a parameter change, we may assume that $s$ is the arc-length parameter of $\gamma$.

\begin{lemma}\label{lambda-const2}
Let $\gamma:I \to \R^4$ be non-degenerate with the arc-length parameter.
Under the notations in Definition \ref{regular-Mannheim}, if $\gamma$ and $\overline{\gamma}$ are second mates or third mates, then $\lambda$ is a non-zero constant.
\end{lemma}
\demo
By differentiating $\overline{\gamma}(s)=\gamma(s)+\lambda(s)\bn_1(s)$,
we have
$$
|\dot{\overline{\gamma}}(s)|\overline{\bt}(s)=(1-\lambda(s)\kappa_1(s))\bt(s)+\lambda'(s)\bn_1(s)+\lambda(s)\kappa_2(s)\bn_2(s).
$$
Since $\bn_1(s)=\pm \overline{\bn}_2(s)$ or $\bn_1(s)=\pm \overline{\bn}_3(s)$, we have $\lambda'(s)=0$ for all $s \in I$. Therefore $\lambda$ is a constant. If $\lambda=0$, then $\overline{\gamma}(s)=\gamma(s)$ for all $s \in I$.
Hence, $\lambda$ is a non-zero constant.
\enD
We give a necessary and sufficient condition of the second type of Mannheim curve for a regular space curve.
\begin{theorem}\label{regular-second-Mannheim-equivalent}
Let $\gamma:I \to \R^4$ be non-degenerate with the arc-length parameter.
Suppose that $\lambda$ is a non-zero constant.
Then $\gamma$ and $\overline{\gamma}$ are second mates with $\overline{\gamma}(s)=\gamma(s)+\lambda \bn_1(s)$ if and only if $\lambda(\kappa_1^2(s)+\kappa_2^2(s))=\kappa_1(s)$,
\begin{eqnarray}
f_1(s) &=& \lambda^2 {\kappa}'_1(s)^2\kappa_2^4(s)\kappa_3^2(s)+(1-\lambda \kappa_1(s))^2{\kappa}'_1(s)^2\kappa_2^2(s)\kappa_3^2(s) \nonumber \\
&&\quad +{\kappa}'_1(s)^2((1-\lambda \kappa_1(s)){\kappa}'_2(s)+\lambda \kappa'_1(s) \kappa_2(s))^2 >0,\\
f_2(s) &=& \kappa_2(s)\left(2{\kappa}'_2(s)^2\kappa_3(s)+\kappa_2(s){\kappa}'_2(s){\kappa}'_3(s)
-\kappa_2(s){\kappa}''_2(s)\kappa_3(s)+\kappa^2_2(s)\kappa_3^3(s)\right) \nonumber\\
&&\quad +\kappa_1(s)(2{\kappa}'_1(s){\kappa}'_2(s)\kappa_3(s)
+{\kappa}'_1(s)\kappa_2(s){\kappa}'_3(s)-{\kappa}''_1(s)\kappa_2(s)\kappa_3(s))=0
\end{eqnarray}
for all $s \in I$.
\end{theorem}
\demo
Suppose that $\overline{\gamma}(s)=\gamma(s)+\lambda \bn_1(s)$ is non-degenerate and  $\bn_1(s)=\pm \overline{\bn}_2(s)$.
Since $\dot{\overline{\gamma}}(s)=(1-\lambda \kappa_1(s))\bt(s)+\lambda \kappa_2(s) \bn_2(s)$,
$|\dot{\overline{\gamma}}(s)|=\sqrt{(1-\lambda \kappa_1(s))^2+\lambda^2 \kappa_2^2(s)}$ and
$$
\overline{\bt}(s)=\frac{1-\lambda \kappa_1(s)}{|\dot{\overline{\gamma}}(s)|}\bt(s)+\frac{\lambda \kappa_2(s)}{|\dot{\overline{\gamma}}(s)|}\bn_2(s).
$$
It follows that
\begin{eqnarray*}
|\dot{\overline{\gamma}}(s)|\overline{\kappa}_1(s)\overline{\bn}_1(s)&=&\dot{\overline{\bt}}(s)\\
&=&\frac{d}{ds}\left(\frac{1-\lambda \kappa_1(s)}{|\dot{\overline{\gamma}}(s)|}\right)\bt(s)+\frac{1}{|\dot{\overline{\gamma}}(s)|}\left( (1-\lambda \kappa_1(s))\kappa_1(s)-\lambda \kappa_2^2(s)\right)\bn_1(s)\\
&&+\frac{d}{ds}\left(\frac{\lambda \kappa_2(s)}{|\dot{\overline{\gamma}}(s)|}\right)\bn_2(s)+\frac{\lambda \kappa_2(s)\kappa_3(s)}{|\dot{\overline{\gamma}}(s)|}\bn_3(s).
\end{eqnarray*}
Since $\bn_1(s)=\pm \overline{\bn}_2(s)$, we have $(1-\lambda \kappa_1(s))\kappa_1(s)-\lambda \kappa_2^2(s)=0$.
Thus, $\lambda (\kappa_1^2(s)+\kappa_2^2(s))=\kappa_1(s)$ for all $s \in I$.
Then we have $(1-\lambda \kappa_1(s))^2+\lambda^2 \kappa_2^2(s)=1-\lambda \kappa_1(s)$.
Therefore, $1-\lambda \kappa_1(s)>0, |\dot{\overline{\gamma}}(s)|=\sqrt{1-\lambda \kappa_1(s)}$ and $\lambda$ is a positive constant.
By a direct calculation, we have
\begin{eqnarray*}
\ddot{\overline{\gamma}}(s)&=&-\lambda {\kappa}'_1(s) \bt(s)+\lambda {\kappa}'_2(s) \bn_2(s)+\lambda \kappa_2(s) \kappa_3(s) \bn_3(s), \\
\dddot{\overline{\gamma}}(s)&=&-\lambda {\kappa}''_1(s) \bt(s)-\frac{{\kappa}'_1(s)}{2} \bn_1(s)+\lambda({\kappa}''_2(s)-\kappa_2(s)\kappa_3^2(s))\bn_2(s)\\
&&+\lambda(2{\kappa}'_2(s)\kappa_3(s)+\kappa_2(s){\kappa}'_3(s))\bn_3(s).
\end{eqnarray*}
Moreover,
\begin{eqnarray*}
&&\dot{\overline{\gamma}}(s) \times \ddot{\overline{\gamma}}(s) \times \dddot{\overline{\gamma}}(s)=\frac{\lambda^2}{2} {\kappa}'_1(s)\kappa_2^2(s) \kappa_3(s) \bt(s)\\
&&\quad +\frac{\lambda^3 \kappa_2(s)}{\kappa_1(s)} \Bigl(\kappa_2(s)(2{\kappa}'_2(s)^2\kappa_3(s)+\kappa_2(s){\kappa}'_2(s){\kappa}'_3(s)
-\kappa_2(s){\kappa}''_2(s)\kappa_3(s)+\kappa_2^2(s)\kappa_3^3(s))\\
&&\quad +\kappa_1(s)(2{\kappa}'_1(s){\kappa}'_2(s)\kappa_3(s)
+{\kappa}'_1(s)\kappa_2(s){\kappa}'_3(s)-{\kappa}''_1(s)\kappa_2(s)\kappa_3(s)) \Bigr)\bn_1(s)\\
&&\quad -\frac{\lambda}{2}(1-\lambda \kappa_1(s)){\kappa}'_1(s)\kappa_2(s)\kappa_3(s) \bn_2(s)\\
&&\quad +\frac{\lambda}{2}\left((1-\lambda \kappa_1(s)) {\kappa}'_1(s) {\kappa}'_2(s)+\lambda {\kappa}'_1(s)^2\kappa_2(s)\right)\bn_3(s).
\end{eqnarray*}
We denote $\dot{\overline{\gamma}}(s) \times \ddot{\overline{\gamma}}(s) \times \dddot{\overline{\gamma}}(s)=A(s)\bt(s)+B(s)\bn_1(s)+C(s)\bn_2(s)+D(s)\bn_3(s)$.
Since
$$
\overline{\bn}_2(s)=-\overline{\bn}_3(s) \times \overline{\bt}(s) \times \overline{\bn}_1(s)
$$
is parallel to
\begin{eqnarray*}
&&-(\dot{\overline{\gamma}}(s) \times \ddot{\overline{\gamma}}(s) \times \dddot{\overline{\gamma}}(s)) \times \dot{\overline{\gamma}}(s)  \times \ddot{\overline{\gamma}}(s)= \lambda^2 \kappa^2_2(s)\kappa_3(s)B(s) \bt(s) \\
&&\quad +\lambda \Bigl((1-\lambda\kappa_1(s))(-{\kappa}'_2(s)D(s)+\kappa_2(s)\kappa_3(s)C(s)) \\
&&\quad -\lambda \kappa_2(s)({\kappa}'_1(s)D(s)+\kappa_2(s)\kappa_3(s)A(s))\Bigr)\bn_1(s)\\
&&\quad -\lambda(1-\lambda\kappa_1(s))\kappa_2(s)\kappa_3(s)B(s)\bn_2(s)\\
&&\quad +\lambda ((1-\lambda\kappa_1(s)){\kappa}'_2(s)+\lambda {\kappa}'_1(s)\kappa_2(s))B(s)\bn_3(s)
\end{eqnarray*}
and $\bn_1(s)=\pm \overline{\bn}_2(s)$, we have $\kappa_3(s)B(s)=0$.
If $B(s) \not=0$, then $\kappa_3(s)=0$ and
$\lambda(1-\lambda\kappa_1(s)){\kappa}'_2(s)+\lambda^2{\kappa}'_1(s)\kappa_2(s)=0$.
It is a contradiction with the fact that $\overline{\bn}_2(s)=\pm \bn_1(s)$.
It follows that $B(s)=0$ for all $s \in I$, that is, $f_2(s)=0$ for all $s \in I$.
Under this condition, $\dot{\overline{\gamma}}(s) \times \ddot{\overline{\gamma}}(s) \times \dddot{\overline{\gamma}}(s) \not=0$ if and only if $f_1(s)>0$ for all $s \in I$.
\par
Conversely, suppose that $\lambda(\kappa_1^2(s)+\kappa_2^2(s))=\kappa_1(s), f_1(s)>0$ and $f_2(s)=0$ for all $s \in I$.
By a direct calculation, $\overline{\gamma}(s)=\gamma(s)+\lambda \bn_1(s)$ is non-degenerate and $\bn_1(s)=\pm \overline{\bn}_2(s)$ for all $s \in I$.
It follows that $\gamma$ and $\overline{\gamma}$ are second mates.
\enD

By the proof of Theorem \ref{regular-second-Mannheim-equivalent}, we have the curvatures of $\overline{\gamma}$.
\begin{proposition}\label{regular-second-Mannheim-curvatures}
Under the same assumptions and notations in Theorem \ref{regular-second-Mannheim-equivalent},
suppose that $\gamma$ and $\overline{\gamma}$ are second mates with  $\overline{\gamma}(s)=\gamma(s)+\lambda \bn_1(s)$.
Then the curvatures of $\overline{\gamma}$ are given by
\begin{eqnarray*}
\overline{\kappa}_1(s) = \frac{\lambda \sqrt{g_1(s)}}{(1-\lambda \kappa_1(s))^{\frac{3}{2}}}, \  \overline{\kappa}_2(s) = \frac{\lambda \sqrt{f_1(s)}}{2(1-\lambda \kappa_1(s))^3 \overline{\kappa}^2_1(s)}, \
\overline{\kappa}_3(s) = \frac{g_2(s)}{(1-\lambda \kappa_1(s))^5 \overline{\kappa}^3_1(s) \overline{\kappa}^2_2(s)},
\end{eqnarray*}
where
{\footnotesize
\begin{eqnarray*}
&&g_1(s)= (1-\lambda \kappa_1(s))({\kappa}'_1(s)^2+{\kappa}'_2(s)^2+\kappa_2^2(s)\kappa^2_3(s))-\frac{{\kappa}'_1(s)^2}{4}, \\
&&g_2(s) = \frac{\lambda^2}{2} {\kappa}'_1(s)\kappa_2^2(s)\kappa_3(s)\left(-\lambda {\kappa}'''_1(s)+\frac{1}{2}\kappa_1(s){\kappa}'_1(s)\right) \\
&& -\frac{1}{2}\left(1-\lambda \kappa_1(s) \right)\lambda {\kappa}'_1(s)\kappa_2(s)\kappa_3(s) \left(-\frac{1}{2}{\kappa}'_1(s)\kappa_2(s)+\lambda {\kappa}'''_2(s)-3\lambda {\kappa}'_2(s)\kappa_3^2(s)-3\lambda \kappa_2(s)\kappa_3(s){\kappa}'_3(s) \right)\\
&& +\frac{\lambda^2}{2} {\kappa}'_1(s)((1-\lambda \kappa_1(s)) {\kappa}'_2(s)+\lambda \kappa_1'(s)\kappa_2(s))(3\lambda {\kappa}''_2(s)\kappa_3(s)-\lambda \kappa_2(s)\kappa_3^3(s)+3\lambda {\kappa}'_2(s) {\kappa}'_3(s)+\lambda \kappa_2(s){\kappa}''_3(s)).
\end{eqnarray*}
}
\end{proposition}
\demo
By Theorem \ref{regular-second-Mannheim-equivalent}, $\lambda (\kappa_1^2(s)+\kappa_2^2(s))=\kappa_1(s)$, $f_1(s)>0$ and $f_2(s)=0$ for all $s \in I$.
By a direct calculation, we have
\begin{eqnarray*}
\dot{\overline{\gamma}}(s)&=&(1-\lambda \kappa_1(s)) \bt(s)+\lambda \kappa_2(s) \bn_2(s), \\
\ddot{\overline{\gamma}}(s)&=&-\lambda {\kappa}'_1(s) \bt(s)+\lambda {\kappa}'_2(s) \bn_2(s)+\lambda \kappa_2(s) \kappa_3(s) \bn_3(s), \\
\dddot{\overline{\gamma}}(s)&=&-\lambda {\kappa}''_1(s) \bt(s)-\frac{{\kappa}'_1(s)}{2} \bn_1(s)+\lambda({\kappa}''_2(s)-\kappa_2(s)\kappa_3^2(s))\bn_2(s)\\
&&+\lambda(2{\kappa}'_2(s)\kappa_3(s)+\kappa_2(s){\kappa}'_3(s))\bn_3(s), \\
{\overline{\gamma}}^{(4)}(s)&=&\left(-\lambda {\kappa}'''_1(s)+\frac{1}{2}\kappa_1(s){\kappa}'_1(s)\right) \bt(s)\\
&&-\left(\lambda {\kappa}''_1(s)\kappa_1(s)+\frac{{\kappa}''_1(s)}{2}+\lambda({\kappa}''_2(s)-\kappa_2(s)\kappa_3^2(s))\kappa_2(s)\right)\bn_1(s)\\
&&+\left(-\frac{1}{2}{\kappa}'_1(s)\kappa_2(s)+\lambda({\kappa}'''_2(s)-3{\kappa}'_2(s)\kappa_3^2(s)- 3\kappa_2(s)\kappa_3(s){\kappa}'_3(s))\right)\bn_2(s)\\
&&+\lambda(3{\kappa}''_2(s)\kappa_3(s)-\kappa_2(s)\kappa_3^3(s)
+3{\kappa}'_2(s){\kappa}'_3(s)+\kappa_2(s){\kappa}''_3(s))\bn_3(s).
\end{eqnarray*}
Therefore, we have $|\dot{\overline{\gamma}}(s)|=\sqrt{1-\lambda\kappa_1(s)}$,
$|\ddot{\overline{\gamma}}(s)|=\lambda \sqrt{{\kappa}'_1(s)^2+{\kappa}'_2(s)^2+{\kappa}_2^2(s)\kappa_3^2(s)}$.
By differentiating $\lambda (\kappa_1^2(s)+\kappa_2^2(s))=\kappa_1(s)$, we have
$$
2\lambda(\kappa_1(s){\kappa}'_1(s)+\kappa_2(s){\kappa}'_2(s))={\kappa}'_1(s).
$$
Since\begin{eqnarray*}
\dot{\overline{\gamma}}(s) \cdot \ddot{\overline{\gamma}}(s) &=&(1-\lambda \kappa_1(s))(-\lambda {\kappa}'_1(s))+\lambda^2 \kappa_2(s){\kappa}'_2(s)=
-\frac{\lambda}{2}{\kappa}'_1(s),
\end{eqnarray*}
we have
$$
g_1(s)=|\dot{\overline{\gamma}}(s)|^2|\ddot{\overline{\gamma}}(s)|^2-(\dot{\overline{\gamma}}(s) \cdot \ddot{\overline{\gamma}}(s))^2= \lambda^2(1-\lambda \kappa_1(s))({\kappa}'_1(s)^2+{\kappa}'_2(s)^2+\kappa_2^2(s)\kappa^2_3(s))-\frac{\lambda^2}{4}{\kappa}'_1(s)^2.
$$
Moreover,
\begin{eqnarray*}
\dot{\overline{\gamma}}(s) \times \ddot{\overline{\gamma}}(s) \times \dddot{\overline{\gamma}}(s)
&=&\frac{\lambda^2}{2} {\kappa}'_1(s)\kappa_2^2(s) \kappa_3(s) \bt(s)\\
&&-\frac{\lambda}{2}(1-\lambda \kappa_1(s)){\kappa}'_1(s)\kappa_2(s)\kappa_3(s) \bn_2(s)\\
&&+\frac{\lambda}{2}\left((1-\lambda \kappa_1(s)) {\kappa}'_1(s) {\kappa}'_2(s)+\lambda {\kappa}'_1(s)^2\kappa_2(s)\right)\bn_3(s),
\end{eqnarray*}
$|\dot{\overline{\gamma}}(s) \times \ddot{\overline{\gamma}}(s) \times \dddot{\overline{\gamma}}(s)|=\lambda \sqrt{f_1(s)}/2$
and
$$
{\rm det}(\dot{\overline{\gamma}}(s),\ddot{\overline{\gamma}}(s),\dddot{\overline{\gamma}}(s),{\overline{\gamma}}^{(4)}(s))=(\dot{\overline{\gamma}}(s) \times \ddot{\overline{\gamma}}(s) \times \dddot{\overline{\gamma}}(s)) \cdot {\overline{\gamma}}^{(4)}(s)=g_2(s).
$$
It follows that
\begin{eqnarray*}
\overline{\kappa}_1(s) &=& \frac{\sqrt{|\dot{\overline{\gamma}}(s)|^2|\ddot{\overline{\gamma}}(s)|^2-(\dot{\overline{\gamma}}(s) \cdot \ddot{\overline{\gamma}}(s))^2}}{|\dot{\overline{\gamma}}(s)|^3}=\frac{\lambda \sqrt{g_1(s)}}{(1-\lambda \kappa_1(s))^\frac{3}{2}}, \\
\overline{\kappa}_2(s) &=& \frac{|\dot{\overline{\gamma}}(s) \times \ddot{\overline{\gamma}}(s) \times \dddot{\overline{\gamma}}(s)|}{|\dot{\overline{\gamma}}(s)|^6 \overline{\kappa}^2_1(s)}=\frac{\lambda \sqrt{f_1(s)}}{2(1-\lambda \kappa_1(s))^3 \overline{\kappa}^2_1(s)}, \\
\overline{\kappa}_3(s) &=& \frac{{\rm det}(\dot{\overline{\gamma}}(s),\ddot{\overline{\gamma}}(s),\dddot{\overline{\gamma}}(s), {\overline{\gamma}}^{(4)}(s) )}{|\dot{\overline{\gamma}}(s)|^{10} \overline{\kappa}^3_1(s) \overline{\kappa}_2^2(s)}=\frac{g_2(s)}{(1-\lambda \kappa_1(s))^5\overline{\kappa}_1^3(s)\overline{\kappa}_2^2(s)}.
\end{eqnarray*}
\enD
We give a necessary and sufficient condition of the third type of Mannheim curve for a regular space curve.
\begin{theorem}\label{regular-third-Mannheim-equivalent}
Let $\gamma:I \to \R^4$ be non-degenerate with the arc-length parameter.
Suppose that $\lambda$ is a non-zero constant.
Then $\gamma$ and $\overline{\gamma}$ are third mates with $\overline{\gamma}(s)=\gamma(s)+\lambda \bn_1(s)$ if and only if $\kappa_1(s)$ and $\kappa_2(s)$ are positive constants with $\lambda(\kappa_1^2+\kappa_2^2)=\kappa_1$ and $\kappa_3(s) \not=0$ for all $s \in I$.
\end{theorem}
\demo
Suppose that $\overline{\gamma}(s)=\gamma(s)+\lambda \bn_1(s)$ is non-degenerate and  $\bn_1(s)=\pm \overline{\bn}_3(s)$.
Since $\dot{\overline{\gamma}}(s)=(1-\lambda \kappa_1(s))\bt(s)+\lambda \kappa_2(s) \bn_2(s)$,
$|\dot{\overline{\gamma}}(s)|=\sqrt{(1-\lambda \kappa_1(s))^2+\lambda^2 \kappa_2^2(s)}$ and
$$
\overline{\bt}(s)=\frac{1-\lambda \kappa_1(s)}{|\dot{\overline{\gamma}}(s)|}\bt(s)+\frac{\lambda \kappa_2(s)}{|\dot{\overline{\gamma}}(s)|}\bn_2(s).
$$
It follows that
\begin{eqnarray*}
|\dot{\overline{\gamma}}(s)|\overline{\kappa}_1(s)\overline{\bn}_1(s)&=&\dot{\overline{\bt}}(s)\\
&=&\frac{d}{ds}\left(\frac{1-\lambda \kappa_1(s)}{|\dot{\overline{\gamma}}(s)|}\right)\bt(s)+\frac{1}{|\dot{\overline{\gamma}}(s)|}\left( (1-\lambda \kappa_1(s))\kappa_1(s)-\lambda \kappa_2^2(s)\right)\bn_1(s)\\
&&+\frac{d}{ds}\left(\frac{\lambda \kappa_2(s)}{|\dot{\overline{\gamma}}(s)|}\right)\bn_2(s)+\frac{\lambda \kappa_2(s)\kappa_3(s)}{|\dot{\overline{\gamma}}(s)|}\bn_3(s).
\end{eqnarray*}
Since $\bn_1(s)=\pm\overline{\bn}_3(s)$, we have $(1-\lambda \kappa_1(s))\kappa_1(s)-\lambda \kappa_2^2(s)=0$.
Thus, $\lambda (\kappa_1^2(s)+\kappa_2^2(s))=\kappa_1(s)$ for all $s \in I$.
Then we have $(1-\lambda \kappa_1(s))^2+\lambda^2 \kappa_2^2(s)=1-\lambda \kappa_1(s)$.
Therefore, $1-\lambda \kappa_1(s)>0$ and $|\dot{\overline{\gamma}}(s)|=\sqrt{1-\lambda \kappa_1(s)}$.
By a direct calculation, we have
\begin{eqnarray*}
\overline{\bt}(s) &=& |\dot{\overline{\gamma}}(s)|\bt(s)+\frac{\lambda \kappa_2(s)}{|\dot{\overline{\gamma}}(s)|}\bn_2(s),\\
\overline{\bn}_1(s) &=& \frac{1}{|\dot{\overline{\gamma}}(s)|\overline{\kappa}_1(s)}\Bigl(
-\frac{\lambda {\kappa}'_1(s)}{2|\dot{\overline{\gamma}}(s)|}\bt(s)+\left(\frac{\lambda {\kappa}'_2(s)}{|\dot{\overline{\gamma}}(s)|}+\frac{\lambda^2 {\kappa}'_1(s)\kappa_2(s)}{2|\dot{\overline{\gamma}}(s)|^3}\right)\bn_2(s)\\
&&+\frac{\lambda \kappa_2(s)\kappa_3(s)}{|\dot{\overline{\gamma}}(s)|}\bn_3(s)
\Bigr).
\end{eqnarray*}
By differentiating $\overline{\bn}_3(s)=\pm \bn_1(s)$, we have $-|\dot{\overline{\gamma}}(s)|\overline{\kappa}_3(s)\overline{\bn}_2(s)=\pm(-\kappa_1(s)\bt(s)+\kappa_2(s)\bn_2(s))$.
It follows that $\overline{\kappa}_3(s) \not=0$ for all $s \in I$ and
\begin{eqnarray*}
\overline{\bn}_2(s)=\mp \frac{1}{|\dot{\overline{\gamma}}(s)|\overline{\kappa}_3(s)}(-\kappa_1(s)\bt(s)+\kappa_2(s)\bn_2(s)).
\end{eqnarray*}
Since
$$
\overline{\bn}_1(s)=\overline{\bn}_2(s) \times \overline{\bn}_3(s) \times \overline{\bt}(s)=
\pm \frac{\kappa_2(s)}{|\dot{\overline{\gamma}}(s)|^2 \overline{\kappa}_3(s)} \bn_3(s),
$$
${\kappa}'_1(s)={\kappa}'_2(s)=0$ for all $s \in I$.
Hence $\kappa_1$ and $\kappa_2$ are positive constants.
It follows that $\lambda$ is also a positive constant.
By a direct calculation, we have
\begin{eqnarray*}
\dot{\overline{\gamma}}(s)&=&(1-\lambda \kappa_1)\bt(s)+\lambda \kappa_2 \bn_2(s), \\
\ddot{\overline{\gamma}}(s)&=&\lambda \kappa_2 \kappa_3(s) \bn_3(s), \\
\dddot{\overline{\gamma}}(s)&=&-\lambda \kappa_2 \kappa_3^2(s) \bn_2(s)+\lambda \kappa_2 {\kappa}'_3(s) \bn_3(s).
\end{eqnarray*}
Therefore, we have
$|\dot{\overline{\gamma}}(s) \times \ddot{\overline{\gamma}}(s) \times \dddot{\overline{\gamma}}(s)|=(1-\lambda \kappa_1)\lambda^2 \kappa_2^2 \kappa_3^2(s) |\kappa_3(s)|$.
By the non-degenerate condition of $\overline{\gamma}$, we have $\kappa_3(s) \not=0$ for all $s \in I$.
\par
Conversely, suppose that $\kappa_1$ and $\kappa_2$ are positive constants with $\lambda(\kappa_1^2+\kappa_2^2)=\kappa_1$ and $\kappa_3(s) \not=0$ for all $s \in I$.
By a direct calculation, $\overline{\gamma}(s)=\gamma(s)+\lambda \bn_1(s)$ is non-degenerate and $\bn_1(s)=\pm \overline{\bn}_3(s)$ for all $s \in I$.
It follows that $\gamma$ and $\overline{\gamma}$ are third mates.
\enD

By the proof of Theorem \ref{regular-third-Mannheim-equivalent}, we have the curvatures of $\overline{\gamma}$.
\begin{proposition}\label{regular-third-Mannheim-curvatures}
Under the same assumptions in Theorem \ref{regular-third-Mannheim-equivalent},
suppose that $\gamma$ and $\overline{\gamma}$ are third mates with  $\overline{\gamma}(s)=\gamma(s)+\lambda \bn_1(s)$.
Then the curvatures of $\overline{\gamma}$ are given by
\begin{eqnarray*}
\overline{\kappa}_1(s) = \frac{\lambda \kappa_2 |\kappa_3(s)|}{1-\lambda \kappa_1}, \
\overline{\kappa}_2(s) = |\kappa_3(s)|, \
\overline{\kappa}_3(s) = \frac{\kappa_2 \kappa_3(s)}{(1-\lambda \kappa_1)|\kappa_3(s)|}.
\end{eqnarray*}
\end{proposition}
\demo
By Theorem \ref{regular-third-Mannheim-equivalent}, $\kappa_1(s)$ and $\kappa_2(s)$ are constants with $\lambda(\kappa_1^2+\kappa_2^2)=\kappa_1$.
By a direct calculation, we have
\begin{eqnarray*}
\dot{\overline{\gamma}}(s)&=&(1-\lambda \kappa_1)\bt(s)+\lambda \kappa_2 \bn_2(s), \\
\ddot{\overline{\gamma}}(s)&=&\lambda \kappa_2 \kappa_3(s) \bn_3(s), \\
\dddot{\overline{\gamma}}(s)&=&-\lambda \kappa_2 \kappa_3^2(s) \bn_2(s)+\lambda \kappa_2 {\kappa}'_3(s) \bn_3(s), \\
{\overline{\gamma}}^{(4)}(s)&=&\lambda \kappa_2^2 {\kappa}_3^2(s) \bn_1(s)-3\lambda \kappa_2\kappa_3(s){\kappa}'_3(s) \bn_2(s) +\lambda (\kappa_2 {\kappa}''_3(s)-\kappa_2 \kappa_3^3(s))\bn_3(s).
\end{eqnarray*}
Therefore, we have
$|\dot{\overline{\gamma}}(s)|=\sqrt{1-\lambda \kappa_1}$,
$|\ddot{\overline{\gamma}}(s)|=\lambda \kappa_2 |\kappa_3(s)|$,
$\dot{\overline{\gamma}}(s) \cdot \ddot{\overline{\gamma}}(s)=0$,
$|\dot{\overline{\gamma}}(s) \times \ddot{\overline{\gamma}}(s) \times \dddot{\overline{\gamma}}(s)|=(1-\lambda \kappa_1)\lambda^2 \kappa_2^2 \kappa_3^2(s) |\kappa_3(s)|$
and
${\rm det}(\dot{\overline{\gamma}}(s),\ddot{\overline{\gamma}}(s),\dddot{\overline{\gamma}}(s),{\overline{\gamma}}^{(4)}(s))=(1-\lambda \kappa_1)\lambda^3 \kappa_2^4\kappa_3^5(s)$.
It follows that
\begin{eqnarray*}
\overline{\kappa}_1(s) &=& \frac{|\ddot{\overline{\gamma}}(s)|}{|\dot{\overline{\gamma}}(s)|^2}=\frac{\lambda \kappa_2 |\kappa_3(s)|}{1-\lambda \kappa_1}, \\
\overline{\kappa}_2(s) &=& \frac{|\dot{\overline{\gamma}}(s) \times \ddot{\overline{\gamma}}(s) \times \dddot{\overline{\gamma}}(s)|}{|\dot{\overline{\gamma}}(s)|^6 \overline{\kappa}^2_1(s)}= |\kappa_3(s)|, \\
\overline{\kappa}_3(s) &=& \frac{{\rm det}(\dot{\overline{\gamma}}(s),\ddot{\overline{\gamma}}(s),\dddot{\overline{\gamma}}(s), {\overline{\gamma}}^{(4)}(s) )}{|\dot{\overline{\gamma}}(s)|^{10} \overline{\kappa}^3_1(s) \overline{\kappa}_2^2(s)}=\frac{\kappa_2 \kappa_3(s)}{(1-\lambda \kappa_1)|\kappa_3(s)|}.
\end{eqnarray*}
\enD
Note that $\overline{\kappa}_3(s)=\pm \kappa_2/(1-\lambda \kappa_1)$ is a constant.

\begin{remark}{\rm
By definitions of second and third types of Mannheim curves,
there is nothing that $\gamma$ is a second and third types of Mannheim curve.
It is also follows from Theorems \ref{regular-second-Mannheim-equivalent} and \ref{regular-third-Mannheim-equivalent}, if $\kappa_1$ and $\kappa_2$ are constants, then $f_1(s)=0$ for all $s \in I$.
}
\end{remark}

\subsection{Framed curves in $\R^4 \times \Delta$}\label{framed-curve}

A framed curve in $\R^4 \times \Delta$ is a smooth space curve with a moving frame, in detail see \cite{Honda-Takahashi1}.
\begin{definition}\label{framed.curve}{\rm
We say that $(\gamma,\nu_1,\nu_2,\nu_3):I \rightarrow \mathbb{R}^4 \times \Delta$ is a {\it framed curve} if $\dot{\gamma}(t) \cdot \nu_i(t)=0$, $i=1,2,3$ for all $t \in I$.
We say that $\gamma:I \to \R^4$ is a {\it framed base curve} if there exists $(\nu_1,\nu_2,\nu_3):I \to \Delta$ such that $(\gamma,\nu_1,\nu_2,\nu_3)$ is a framed curve.
}
\end{definition}

We denote $\bmu(t) = \nu_1(t) \times \nu_2(t) \times \nu_3(t)$.
Then $\{ \nu_1(t),\nu_2(t), \nu_3(t), \bmu(t) \}$ is a moving frame along the framed base curve $\gamma(t)$ in $\R^4$ and we have the Frenet-Serret type formula:
$$
\left(
\begin{array}{c}
\dot{\nu_1}(t)\\
\dot{\nu_2}(t)\\
\dot{\nu_3}(t)\\
\dot{\bmu}(t)
\end{array} \right)=
\left(
\begin{array}{cccc}
0 & \ell_1(t) & \ell_2(t) & \ell_3(t)\\
-\ell_1(t) & 0 & \ell_4(t) & \ell_5(t)\\
-\ell_2(t) & -\ell_4(t) & 0 & \ell_6(t)\\
-\ell_3(t) & -\ell_5(t) & -\ell_6(t) &0
\end{array}\right)
\left(
\begin{array}{c}
\nu_1(t)\\
\nu_2(t)\\
\nu_3(t)\\
\bmu(t)
\end{array}\right),
\ \dot{\gamma}(t)=\alpha(t)\bmu(t),
$$
where
\begin{eqnarray*}
&& \ell_1(t) = \dot{\nu_1}(t) \cdot \nu_2(t), \ \ell_2(t) = \dot{\nu_1}(t) \cdot \nu_3(t), \ \ell_3(t) = \dot{\nu_1}(t) \cdot \bmu(t), \\
&& \ell_4(t) = \dot{\nu_2}(t) \cdot \nu_3(t), \ \ell_5(t) = \dot{\nu_2}(t) \cdot \bmu(t), \ \ell_6(t) = \dot{\nu_3}(t) \cdot \bmu(t)
\end{eqnarray*}
and $\alpha(t)=\dot{\gamma}(t) \cdot \bmu(t)$.
We call the mapping $(\ell_1,\dots,\ell_6,\alpha): I \to \R^7$ {\it the curvature of the framed curve} $(\gamma,\nu_1,\nu_2,\nu_3): I \to \R^4 \times \Delta$.
Note that $t_0$ is a singular point of $\gamma$ if and only if $\alpha(t_0) = 0$.
\begin{definition}
{\rm
Let $(\gamma,\nu_1,\nu_2,\nu_3)$ and $(\widetilde{\gamma},\widetilde{\nu}_1,\widetilde{\nu}_2,\widetilde{\nu}_3):I \rightarrow \mathbb{R}^4 \times \Delta$ be framed curves.
We say that $(\gamma,\nu_1,\nu_2,\nu_3)$ and $(\widetilde{\gamma},\widetilde{\nu}_1,\widetilde{\nu}_2,\widetilde{\nu}_3)$ are {\it congruent as framed curves} if there exist a constant rotation $A \in SO(4)$ and a translation $\ba \in \mathbb{R}^4$ such that $\widetilde{\gamma}(t) = A(\gamma(t)) +\ba$ and $\widetilde{\nu}_i(t) = A(\nu_i(t))$,  $i=1,2,3$ for all $t \in I$.
}
\end{definition}
\par
We gave the existence and uniqueness theorems for framed curves in terms of the curvatures in \cite{Honda-Takahashi1}.
\begin{theorem}[Existence Theorem for framed curves \cite{Honda-Takahashi1}]\label{existence.framed}
Let $(\ell_1,\dots,\ell_6,\alpha):I \rightarrow \mathbb{R}^7$ be a smooth mapping.
Then there exists a framed curve $(\gamma,\nu_1,\nu_2,\nu_3):I \to \R^4 \times \Delta$ whose curvature is given by $(\ell_1,\dots,\ell_6,\alpha)$.
\end{theorem}
\begin{theorem}[Uniqueness Theorem for framed curves \cite{Honda-Takahashi1}]\label{uniqueness.framed}
Let $(\gamma,\nu_1,\nu_2,\nu_3)$ and \\ $(\widetilde{\gamma},\widetilde{\nu}_1,\widetilde{\nu}_2,\widetilde{\nu}_3):I \to \R^4 \times \Delta$ be framed curves with curvatures $(\ell_1,\dots,\ell_6,\alpha)$ and $(\widetilde{\ell}_1,\dots,\widetilde{\ell}_6,\widetilde{\alpha})$, respectively.
Then $(\gamma,\nu_1,\nu_2,\nu_3)$ and $(\widetilde{\gamma},\widetilde{\nu}_1,\widetilde{\nu}_2,\widetilde{\nu}_3)$ are congruent as framed curves if and only if the curvatures $(\ell_1,\dots,\ell_6, \alpha)$ and $(\widetilde{\ell}_1, \dots, \widetilde{\ell}_6, \widetilde{\alpha})$ coincide.
\end{theorem}

Let $(\gamma,\nu_1,\nu_2,\nu_3): I \to \R^4 \times \Delta$ be a framed curve with the curvature $(\ell_1,\dots,\ell_6,\alpha)$.
For the normal hyperplane of $\gamma(t)$, spanned by $\nu_1(t), \nu_2(t), \nu_3(t)$, there is some ambient of framed curves similarly to the case of the Bishop frame of a regular space curve (cf. \cite{Bishop}).
We define $(\widetilde{\nu}_1(t), \widetilde{\nu}_2(t), \widetilde{\nu}_3(t)) \in \Delta$ by $A(t) \in SO(3)$,
\begin{eqnarray*}\label{eqn-rotation}
\left(
\begin{array}{c}
\widetilde{\nu}_1(t)\\
\widetilde{\nu}_2(t)\\
\widetilde{\nu}_3(t)
\end{array}
\right)
=A(t)
\left(
\begin{array}{c}
{\nu}_1(t)\\
\nu_2(t)\\
\nu_3(t)
\end{array}
\right),
\end{eqnarray*}
where $A(t)=$
{\footnotesize
$$
\left(
\begin{array}{ccc}
\cos \phi(t) \cos \psi(t)- \sin \phi(t)\cos \theta(t)\sin \psi(t)&-\cos \phi(t) \sin \psi(t)-\sin \phi(t)\cos \theta(t)\cos \psi(t)&\sin \phi(t) \sin \theta(t)\\
\sin \phi(t) \cos \psi(t)+ \cos \phi(t)\cos \theta(t)\sin \psi(t)&-\sin \phi(t) \sin \psi(t)+\cos \phi(t) \cos \theta(t) \cos \psi(t)&-\cos \phi(t) \sin \theta(t)\\
\sin \theta(t)\sin \psi(t) & \sin \theta(t)\cos \psi(t) &\cos \theta(t)
\end{array}
\right)
$$}
and $\phi, \psi, \theta:I \to \R$ are smooth functions (cf. \cite{Gelfand-Minlos-Shapiro}).
Then $(\gamma, \widetilde{\nu}_1, \widetilde{\nu}_2, \widetilde{\nu}_3) : I \rightarrow \R^4 \times \Delta $ is also a framed curve and $\widetilde{\bmu}(t) = \bmu(t)$. By a direct calculation, we have
the Frenet-Serret type formula:
$$
\left(
\begin{array}{c}
\dot{\widetilde{\nu}_1}(t)\\
\dot{\widetilde{\nu}_2}(t)\\
\dot{\widetilde{\nu}_3}(t)\\
\dot{{\bmu}}(t)
\end{array} \right)=
\left(
\begin{array}{cccc}
0 & \widetilde{\ell}_1(t) & \widetilde{\ell}_2(t) & \widetilde{\ell}_3(t)\\
-\widetilde{\ell}_1(t) & 0 & \widetilde{\ell}_4(t) & \widetilde{\ell}_5(t)\\
-\widetilde{\ell}_2(t) & -\widetilde{\ell}_4(t) & 0 & \widetilde{\ell}_6(t)\\
-\widetilde{\ell}_3(t) & -\widetilde{\ell}_5(t) & -\widetilde{\ell}_6(t) &0
\end{array}\right)
\left(
\begin{array}{c}
\widetilde{\nu}_1(t)\\
\widetilde{\nu}_2(t)\\
\widetilde{\nu}_3(t)\\
{\bmu}(t)
\end{array}\right),
\ \dot{\gamma}(t)=\alpha(t){\bmu}(t),
$$
where
\begin{eqnarray*}
\widetilde{\ell}_1(t)&=&-\dot{\phi}(t)+(\ell_1(t)-\dot{\psi}(t))\cos\theta(t)+(-\ell_2(t)\cos\psi(t)+\ell_4(t)\sin\psi(t))\sin\theta(t),\label{ell_11}\\
\widetilde{\ell}_2(t)&=&\left(\dot{\theta}(t)-(\ell_2(t)\sin\psi(t)+\ell_4(t)\cos\psi(t))\right)\sin\phi(t) \nonumber\\
&&+\left((\ell_2(t)\cos\psi(t)-\ell_4(t)\sin\psi(t))\cos\theta(t)+(\ell_1(t)-\dot{\psi}(t))\sin\theta(t)\right)\cos\phi(t),\label{ell_22}\\
\widetilde{\ell}_3(t)&=&\left(\ell_6(t)\sin\theta(t)-(\ell_3(t)\sin\psi(t)+\ell_5(t)\cos\psi(t))\cos\theta(t)\right)\sin\phi(t)\nonumber\\
&&+(\ell_3(t)\cos\psi(t)-\ell_5(t)\sin\psi(t))\cos\phi(t),\label{ell_33}\\
\widetilde{\ell}_4(t)&=&\left(-\dot{\theta}(t)+(\ell_2(t)\sin\psi(t)+\ell_4(t)\cos\psi(t))\right)\cos\phi(t)\nonumber\\
&&+\left((\ell_1(t)-\dot{\psi}(t))\sin\theta(t)+(\ell_2(t)\cos\psi(t)-\ell_4(t)\sin\psi(t))\cos\theta(t)\right)\sin\phi(t),\label{ell_44}\\
\widetilde{\ell}_5(t)&=&(-\ell_6(t)\sin\theta(t)+(\ell_3(t)\sin\psi(t)+\ell_5(t)\cos\psi(t))\cos\theta(t))\cos\phi(t)\nonumber\\
&&+(\ell_3(t)\cos\psi(t)-\ell_5(t)\sin\psi(t))\sin\phi(t),\label{ell_55}\\
\widetilde{\ell}_6(t)&=&(\ell_3(t)\sin\psi(t)+\ell_5(t)\cos\psi(t))\sin\theta(t)+\ell_6(t)\cos\theta(t).\label{ell_66}
\end{eqnarray*}

\begin{corollary}
Under the above notations,
if we take smooth functions $\phi, \psi, \theta:I \to \R$ which satisfy $\sin\theta(t)\neq0$ for all $t\in I$ and
\begin{eqnarray*}
\dot{\phi}(t)&=&(-\ell_2(t)\cos\psi(t)+\ell_4(t)\sin\psi(t))/\sin \theta(t),\\
\dot{\psi}(t)&=&(\ell_2(t)\cos\psi(t)-\ell_4(t)\sin\psi(t))\cos \theta(t)/\sin \theta(t)+\ell_1(t),\\
\dot{\theta}(t)&=&\ell_2(t)\sin\psi(t)+\ell_4(t)\cos\psi(t),
\end{eqnarray*}
then $\widetilde{\ell}_1(t)=\widetilde{\ell}_2(t)=\widetilde{\ell}_4(t)=0$ for all $t \in I$.
\end{corollary}

We call the frame $\{\widetilde{\nu}_1(t),\widetilde{\nu}_2(t),\widetilde{\nu}_3(t),{\bmu}(t)\}$ an {\it adapted frame} along $\gamma(t)$.
In this case, the Frenet-Serret type formula is given by
\begin{eqnarray}\label{adapt-frame}
\left(
\begin{array}{c}
\dot{\widetilde{\nu}_1}(t)\\
\dot{\widetilde{\nu}_2}(t)\\
\dot{\widetilde{\nu}_3}(t)\\
\dot{{\bmu}}(t)
\end{array} \right)=
\left(
\begin{array}{cccc}
0 & 0 & 0 & \widetilde{\ell}_3(t)\\
0 & 0 & 0 & \widetilde{\ell}_5(t)\\
0 & 0 & 0 & \widetilde{\ell}_6(t)\\
-\widetilde{\ell}_3(t) & -\widetilde{\ell}_5(t) & -\widetilde{\ell}_6(t) &0
\end{array}\right)
\left(
\begin{array}{c}
\widetilde{\nu}_1(t)\\
\widetilde{\nu}_2(t)\\
\widetilde{\nu}_3(t)\\
{\bmu}(t)
\end{array}\right), \
\dot{\gamma}(t)=\alpha(t)\bmu(t).
\end{eqnarray}
See also \cite{YKT}.

\section{Bertrand curves of framed curves}\label{Bertrand-framed}

Let $(\gamma,\nu_1,\nu_2,\nu_3)$ and $(\overline{\gamma},\overline{\nu}_1,\overline{\nu}_2,\overline{\nu}_3):I \to \R^4 \times \Delta$ be framed curves with the curvatures $(\ell_1,\dots,\ell_6,\alpha)$ and $(\overline{\ell}_1,\dots,\overline{\ell}_6,\overline{\alpha})$, respectively.
Suppose that $\gamma$ and $\overline{\gamma}$ are different curves, that is, $\gamma \not\equiv \overline{\gamma}$.

\begin{definition}\label{Bertrand.framed}{\rm
We say that framed curves $(\gamma,\nu_1,\nu_2,\nu_3)$ and $(\overline{\gamma},\overline{\nu}_1,\overline{\nu}_2,\overline{\nu}_3)$ are {\it Bertrand mates} ({\it first mates} or, $(\nu_1,\overline{\nu}_1)$-{\it mates}) if there exists a smooth function $\lambda:I \to \R$ such that $\overline{\gamma}(t)=\gamma(t)+\lambda(t)\nu_1(t)$ and $\nu_1(t)=\overline{\nu}_1(t)$ for all $t \in I$.
We also say that $(\gamma,\nu_1,\nu_2,\nu_3):I \to \R^4 \times \Delta$ is a {\it Bertrand curve} if there exists a framed curve  $(\overline{\gamma},\overline{\nu}_1,\overline{\nu}_2,\overline{\nu}_3):I \to \R^4 \times \Delta$ such that $(\gamma,\nu_1,\nu_2,\nu_3)$ and $(\overline{\gamma},\overline{\nu}_1,\overline{\nu}_2,\overline{\nu}_3)$ are Bertrand mates.
}
\end{definition}

\begin{lemma}\label{lambda-const}
Under the notations in Definition \ref{Bertrand.framed},
if $(\gamma,\nu_1,\nu_2,\nu_3)$ and $(\overline{\gamma},\overline{\nu}_1,\overline{\nu}_2,\overline{\nu}_3)$ are Bertrand mates, then $\lambda$ is a non-zero constant.
\end{lemma}
\demo
By differentiating $\overline{\gamma}(t)=\gamma(t)+\lambda(t)\nu_1(t)$, we have
$$
\overline{\alpha}(t) \overline{\bmu}(t)
=(\alpha(t)+\lambda(t)\ell_3(t))\bmu(t)+\dot{\lambda}(t)\nu_1(t)+\lambda(t)\ell_1(t)\nu_2(t)+\lambda(t)\ell_2(t)\nu_3(t)
$$
for all $t \in I$.
Since $\overline{\nu}_1(t)=\nu_1(t)$, we have $\dot{\lambda}(t)=0$ for all $t \in I$.
Therefore $\lambda$ is a constant.
If $\lambda=0$, then $\overline{\gamma}(t)=\gamma(t)$ for all $t \in I$.
Hence, $\lambda$ is a non-zero constant.
\enD

We give a necessary and sufficient condition of the Bertrand curve for a framed curve.

\begin{theorem}\label{framed-Bertrand-equivalent}
Let $(\gamma,\nu_1,\nu_2,\nu_3):I \to \R^4 \times \Delta$ be a framed curve with the curvature $(\ell_1,\dots,\ell_6,\alpha)$.
Then $(\gamma,\nu_1,\nu_2,\nu_3)$ is a Bertrand curve if and only if there exist a non-zero  constant $\lambda$ and smooth functions $\psi, \theta:I \to \R$ such that
\begin{eqnarray}
&&\ell_1(t)\cos \psi(t)-\ell_2(t)\sin \psi(t)=0, \label{Bertrand-condition1}\\
&&\lambda(\ell_1(t) \sin \psi(t)+ \ell_2(t)\cos \psi(t))\cos \theta(t)-(\alpha(t)+\lambda \ell_3(t))\sin \theta(t)=0 \label{Bertrand-condition2}
\end{eqnarray}
for all $t \in I$.
\end{theorem}
\demo
Suppose that $(\gamma,\nu_1,\nu_2,\nu_3)$ is a Bertrand curve.
By Definition \ref{Bertrand.framed} and Lemma \ref{lambda-const}, there exist a framed curve $(\overline{\gamma},\overline{\nu}_1,\overline{\nu}_2,\overline{\nu}_3)$ and a non-zero constant $\lambda \in \R$ such that $\overline{\gamma}(t)=\gamma(t)+\lambda \nu_1(t)$ and $\nu_1(t)=\overline{\nu}_1(t)$ for all $t \in I$.
By differentiating $\overline{\gamma}(t)=\gamma(t)+\lambda \nu_1(t)$, we have
$$
\overline{\alpha}(t) \overline{\bmu}(t)=(\alpha(t)+\lambda \ell_3(t))\bmu(t)+\lambda \ell_1(t) \nu_2(t)+\lambda \ell_2(t)\nu_3(t).
$$
Since $\nu_1(t)=\overline{\nu}_1(t)$, there exists $A(t) \in SO(3)$ such that
\begin{eqnarray*}\label{eqn-rotation}
\left(
\begin{array}{c}
\overline{\nu}_2(t)\\
\overline{\nu}_3(t)\\
\overline{\bmu}(t)
\end{array}
\right)
=A(t)
\left(
\begin{array}{c}
{\nu}_2(t)\\
\nu_3(t)\\
{\bmu}(t)
\end{array}
\right),
\end{eqnarray*}
where $A(t)=$
{\footnotesize
$$
\left(
\begin{array}{ccc}
\cos \phi(t) \cos \psi(t)- \sin \phi(t)\cos \theta(t)\sin \psi(t)&-\cos \phi(t) \sin \psi(t)-\sin \phi(t)\cos \theta(t)\cos \psi(t)&\sin \phi(t) \sin \theta(t)\\
\sin \phi(t) \cos \psi(t)+ \cos \phi(t)\cos \theta(t)\sin \psi(t)&-\sin \phi(t) \sin \psi(t)+\cos \phi(t) \cos \theta(t) \cos \psi(t)&-\cos \phi(t) \sin \theta(t)\\
\sin \theta(t)\sin \psi(t) & \sin \theta(t)\cos \psi(t) &\cos \theta(t)
\end{array}
\right)
$$}
and $\phi, \psi, \theta:I \to \R$ are smooth functions.
Then we have
$$
\overline{\bmu}(t)=\sin \theta(t)\sin \psi(t) \nu_2(t)+ \sin \theta(t)\cos \psi(t) \nu_3(t)+\cos \theta(t)\bmu(t).
$$
It follows that
\begin{eqnarray}
\overline{\alpha}(t) \sin \theta(t) \sin \psi(t)&=&\lambda \ell_1(t), \label{ell_1}\\
\overline{\alpha}(t) \sin \theta(t) \cos \psi(t)&=&\lambda \ell_2(t), \label{ell_2}\\
\overline{\alpha}(t) \cos \theta(t) &=& \alpha(t)+\lambda \ell_3(t). \label{ell_3}
\end{eqnarray}
Since \eqref{ell_1} and \eqref{ell_2}, we have
\begin{eqnarray}
&& \ell_1(t) \cos \psi(t)-\ell_2(t)\sin \psi(t)=0, \nonumber\\
&& \overline{\alpha}(t) \sin \theta(t) =\lambda (\ell_1(t) \sin \psi(t)+\ell_2(t)\cos\psi(t)). \label{alpha}
\end{eqnarray}
Since \eqref{ell_3} and \eqref{alpha}, we also have
$$
\lambda(\ell_1(t) \sin \psi(t)+ \ell_2(t)\cos \psi(t))\cos \theta(t)-(\alpha(t)+\lambda \ell_3(t))\sin \theta(t)=0
$$
for all $t \in I$.
\par
Conversely, there exist a non-zero  constant $\lambda$ and smooth functions $\psi, \theta:I \to \R$ such that \eqref{Bertrand-condition1} and \eqref{Bertrand-condition2} hold.
We define a mapping $(\overline{\gamma},\overline{\nu}_1,\overline{\nu}_2,\overline{\nu}_3):I \to \R^4 \times \Delta$ by
\begin{eqnarray}
\overline{\gamma}(t)&=&\gamma(t)+\lambda \nu_1(t), \label{Bertrand-1}\\
\overline{\nu}_1(t)&=&\nu_1(t), \label{Bertrand-2}\\
\overline{\nu}_2(t)&=&(\cos\phi(t)\cos\psi(t)-\sin\phi(t)\cos\theta(t)\sin\psi(t))\nu_2(t)\notag\\
&&+(-\cos\phi(t)\sin\psi(t)-\sin\phi(t)\cos\theta(t)\cos\psi(t))\nu_3(t)+\sin\phi(t)\sin\theta(t)\bmu(t), \label{Bertrand-3}\\
\overline{\nu}_3(t)&=&(\sin\phi(t)\cos\psi(t)+\cos\phi(t)\cos\theta(t)\sin\psi(t))\nu_2(t)\notag\\
&&+(-\sin\phi(t)\sin\psi(t)+\cos\phi(t)\cos\theta(t)\cos\psi(t))\nu_3(t)-\cos\phi(t)\sin\theta(t)\bmu(t), \label{Bertrand-4}
\end{eqnarray}
where $\phi:I \to \R$ is a smooth function.
By a direct calculation, $(\overline{\gamma},\overline{\nu}_1,\overline{\nu}_2,\overline{\nu}_3)$ is a framed curve.
Therefore, $(\gamma,\nu_1,\nu_2,\nu_3)$ is a Bertrand curve.
\enD
\begin{corollary}\label{ell-Bertand}
Let $(\gamma,\nu_1,\nu_2,\nu_3):I \to \R^4 \times \Delta$ be a framed curve with the curvature $(\ell_1,\dots,\ell_6,\alpha)$.
If $\ell_1(t)=\ell_2(t)=0$ for all $t \in I$, then $(\gamma,\nu_1,\nu_2,\nu_3)$ is a Bertrand curve.
\end{corollary}
\demo
If $\ell_1(t)=\ell_2(t)=0$ for all $t \in I$, then equation \eqref{Bertrand-condition1} is satisfied.
If we take $\theta(t)=0$ for all $t \in I$, then equation \eqref{Bertrand-condition2} is satisfied.
By Theorem \ref{framed-Bertrand-equivalent}, $(\gamma,\nu_1,\nu_2,\nu_3)$ is a Bertrand curve.
\enD
Let $(\gamma,\nu_1,\nu_2,\nu_3):I \to \R^4 \times \Delta$ be a framed curve with the curvature $(\ell_1,\dots,\ell_6,\alpha)$.
If we take an adapted frame $\{\widetilde{\nu}_1, \widetilde{\nu}_2, \widetilde{\nu}_3,{\bmu}\}$,
then the curvature is given by $(0,0, \widetilde{\ell}_3, 0, \widetilde{\ell}_5,\widetilde{\ell}_6,\alpha)$, see \eqref{adapt-frame}. By Theorem \ref{framed-Bertrand-equivalent} or Corollary \ref{ell-Bertand}, we have the following.
\begin{corollary}
For an adapted frame, $(\gamma,\widetilde{\nu}_1,\widetilde{\nu}_2,\widetilde{\nu}_3)$ is always a Bertrand curve.
\end{corollary}
\begin{proposition}\label{Bertrand-frame-curvature}
Suppose that $(\gamma,\nu_1,\nu_2,\nu_3)$ and $(\overline{\gamma},\overline{\nu}_1,\overline{\nu}_2,\overline{\nu}_3):I \to \R^4 \times \Delta$ are Bertrand mates which is given by \eqref{Bertrand-1}, \eqref{Bertrand-2}, \eqref{Bertrand-3} and \eqref{Bertrand-4}.
Then the curvature $(\overline{\ell}_1,\dots,\overline{\ell}_6,\overline{\alpha})$ of $(\overline{\gamma},\overline{\nu}_1,\overline{\nu}_2,\overline{\nu}_3)$ is given by
\begin{eqnarray*}\label{Bertrand-curvature}
\overline{\ell}_1(t) &=& -\alpha(t)\sin\phi(t)\sin\theta(t)/\lambda ,\\
\overline{\ell}_2(t) &=& (\ell_1(t)\sin\psi(t)+\ell_2(t)\cos\psi(t))\cos\phi(t) \cos\theta(t)-\ell_3(t) \cos \phi(t)\sin\theta(t),\\
\overline{\ell}_3(t) &=& (\ell_1(t) \sin \psi(t)+\ell_2(t)\cos \psi(t))\sin \theta(t)+\ell_3(t)\cos \theta(t),\\
\overline{\ell}_4(t) &=& -\dot{\phi}(t)+(\ell_4(t)-\dot{\psi}(t))\cos \theta(t)-(\ell_5(t)\cos \psi(t)-\ell_6(t)\sin \psi(t))\sin \theta(t),\\
\overline{\ell}_5(t) &=& \dot{\theta}(t)\sin\phi(t)+(\ell_4(t)-\dot{\psi}(t))\cos\phi(t)\sin \theta(t)-(\ell_5(t)\sin\psi(t)+\ell_6(t)\cos\psi(t))\sin\phi(t)\\
&&+(\ell_5(t)\cos \psi(t)-\ell_6(t)\sin \psi(t))\cos\phi(t)\cos \theta(t),\\
\overline{\ell}_6(t) &=& -\dot{\theta}(t)\cos\phi(t)+(\ell_4(t)-\dot{\psi}(t))\sin\phi(t)\sin\theta(t)
+(\ell_5(t)\sin \psi(t)+\ell_6(t)\cos \psi(t))\cos\phi(t)\\&&+(\ell_5(t)\cos\psi(t)-\ell_6(t)\sin\psi(t))\sin\phi(t)\cos\theta(t),\\
\overline{\alpha}(t) &=& \lambda (\ell_1(t)\sin\psi(t)+\ell_2(t)\cos \psi(t))\sin \theta(t)+(\alpha(t)+\lambda \ell_3(t))\cos \theta(t).
\end{eqnarray*}
\end{proposition}
\demo
By a direct calculation, we have
$$
\overline{\bmu}(t)=\sin \theta(t) \sin \psi(t) \nu_2(t)+\sin \theta(t)\cos \psi(t)\nu_3(t)+\cos \theta(t)\bmu(t).
$$
By differentiating $\overline{\gamma},\overline{\nu}_1,\overline{\nu}_2$ and $\overline{\nu}_3$,
\begin{eqnarray*}
&&\dot{\overline{\gamma}}(t)=\lambda \ell_1(t)\nu_2(t)+\lambda \ell_2(t)\nu_3(t)+(\alpha(t)+\lambda\ell_3(t))\bmu(t),\\
&&\dot{\overline{\nu}}_1(t)=\ell_1(t)\nu_2(t)+\ell_2(t)\nu_3(t)+\ell_3(t)\bmu(t),
\end{eqnarray*}
\begin{eqnarray*}
&&\dot{\overline{\nu}}_2(t)=\big(\ell_1(t)(\sin\phi(t)\cos\theta(t)\sin\psi(t)-\cos\phi(t)\cos\psi(t))\\
&&\quad+\ell_2(t)(\sin\phi(t)\cos\theta(t)\cos\psi(t)+\cos\phi(t)\sin\psi(t))-\ell_3(t)\sin\phi(t)\sin\theta(t)\big)\nu_1(t)\\
&&\quad+\big(-\dot{\phi}(t)\sin\phi(t)\cos\psi(t)-\dot{\psi}(t)\cos\phi(t)\sin\psi(t)+\dot{\theta}(t)\sin\phi(t)\sin \theta(t)\sin \psi(t)\\
&&\quad-\dot{\phi}(t)\cos\phi(t)\cos\theta(t)\sin\psi(t)-\dot{\psi}(t)\sin\phi(t)\cos\theta(t)\cos\psi(t)\\
&&\quad+\ell_4(t)(\cos\phi(t)\sin\psi(t)+\sin\phi(t)\cos\theta(t)\cos \psi(t))-\ell_5(t)\sin\phi(t)\sin\theta(t)\big)\nu_2(t)\\
&&\quad+\big(\dot{\phi}(t)\sin\phi(t)\sin\psi(t)-\dot{\psi}(t)\cos\phi(t)\cos\psi(t)+\dot{\theta}(t)\sin\phi(t)\sin\theta(t)\cos\psi(t)\\
&&\quad-\dot{\phi}(t)\cos\phi(t)\cos\theta(t)\cos\psi(t)+\dot{\psi}(t)\sin\phi(t)\cos\theta(t)\sin\psi(t)\\
&&\quad+\ell_4(t)(\cos\phi(t)\cos\psi(t)-\sin\phi(t)\cos\theta(t)\sin\psi(t))-\ell_6(t)\sin\phi(t)\sin\theta(t)\big)\nu_3(t)\\
&&\quad+\big(\dot{\phi}(t)\cos\phi(t)\sin\theta(t)+\dot{\theta}(t)\sin\phi(t)\cos\theta(t)\\
&&\quad+\ell_5(t)(\cos\phi(t)\cos\psi(t)-\sin\phi(t)\cos\theta(t)\sin\psi(t))\\
&&\quad-\ell_6(t)(\cos\phi(t)\sin\psi(t)+\sin\phi(t)\cos\theta(t)\cos\psi(t))\big)\bmu(t),
\end{eqnarray*}
\begin{eqnarray*}
&&\dot{\overline{\nu}}_3(t)=\big(-\ell_1(t)(\cos\phi(t)\cos\theta(t)\sin\psi(t)+\sin\phi(t)\cos\psi(t))\\
&&\quad+\ell_2(t)(-\cos\phi(t)\cos \theta(t)\cos\psi(t)+\sin\phi(t)\sin\psi(t))+\ell_3(t)\cos\phi(t)\sin \theta(t)\big)\nu_1(t)\\
&&\quad+\big(\dot{\phi}(t)\cos\phi(t)\cos\psi(t)-\dot{\psi}(t)\sin\phi(t)\sin\psi(t)-\dot{\theta}(t)\cos\phi(t)\sin \theta(t)\sin \psi(t)\\
&&\quad-\dot{\phi}(t)\sin\phi(t)\cos\theta(t)\sin\psi(t)+\dot{\psi}(t)\cos\phi(t)\cos\theta(t)\cos\psi(t)\\
&&\quad+\ell_4(t)(\sin\phi(t)\sin\psi(t)-\cos\phi(t)\cos\theta(t)\cos \psi(t))+\ell_5(t)\cos\phi(t)\sin \theta(t)\big)\nu_2(t)\\
&&\quad+\big(-\dot{\phi}(t)\cos\phi(t)\sin\psi(t)-\dot{\psi}(t)\sin\phi(t)\cos\psi(t)-\dot{\theta}(t)\cos\phi(t)\sin\theta(t)\cos\psi(t)\\
&&\quad-\dot{\phi}(t)\sin\phi(t)\cos\theta(t)\cos\psi(t)-\dot{\psi}(t)\cos\phi(t)\cos\theta(t)\sin\psi(t)\\
&&\quad+\ell_4(t)(\sin\phi(t)\cos\psi(t)+\cos\phi(t)\cos\theta(t)\sin \psi(t))+\ell_6(t)\cos\phi(t)\sin \theta(t)\big)\nu_3(t)\\
&&\quad+\big(\dot{\phi}(t)\sin\phi(t)\sin\theta(t)-\dot{\theta}(t)\cos\phi(t)\cos\theta(t)\\
&&\quad+\ell_5(t)(\sin\phi(t)\cos\psi(t)+\cos\phi(t)\cos\theta(t)\sin\psi(t))\\
&&\quad+\ell_6(t)(-\sin\phi(t)\sin\psi(t)+\cos\phi(t)\cos\theta(t)\cos\psi(t))\big)\bmu(t).
\end{eqnarray*}
By a direct calculation, equations \eqref{Bertrand-condition1} and \eqref{Bertrand-condition2}, we have the curvature.
\enD

In particular, if we take $\phi(t)=0$ for all $t\in I$, then $(\overline{\gamma},\overline{\nu}_1,\overline{\nu}_2,\overline{\nu}_3):I \to \R^4 \times \Delta$ is given by
\begin{eqnarray}
\overline{\gamma}(t)&=&\gamma(t)+\lambda \nu_1(t), \label{Bertrand-5}\\
\overline{\nu}_1(t) &=&\nu_1(t), \label{Bertrand-6}\\
\overline{\nu}_2(t) &=& \cos \psi(t) \nu_2(t)-\sin \psi(t)\nu_3(t), \label{Bertrand-7}\\
\overline{\nu}_3(t) &=& \cos \theta(t) \sin \psi(t) \nu_2(t)+\cos \theta(t) \cos \psi(t)\nu_3(t)-\sin \theta(t) \bmu(t). \label{Bertrand-8}
\end{eqnarray}
It is easy to see that $(\overline{\gamma},\overline{\nu}_1,\overline{\nu}_2,\overline{\nu}_3)$ is a framed curve.
Moreover, $(\gamma,\nu_1,\nu_2,\nu_3)$ and $(\overline{\gamma},\overline{\nu}_1,\overline{\nu}_2,\overline{\nu}_3)$ are Bertrand mates.

\begin{corollary}\label{Bertrand-frame-curvature}
Suppose that $(\gamma,\nu_1,\nu_2,\nu_3)$ and $(\overline{\gamma},\overline{\nu}_1,\overline{\nu}_2,\overline{\nu}_3):I \to \R^4 \times \Delta$ are Bertrand mates which is given by \eqref{Bertrand-5}, \eqref{Bertrand-6}, \eqref{Bertrand-7} and \eqref{Bertrand-8}.
Then the curvature $(\overline{\ell}_1,\dots,\overline{\ell}_6,\overline{\alpha})$ of $(\overline{\gamma},\overline{\nu}_1,\overline{\nu}_2,\overline{\nu}_3)$ is given by
\begin{eqnarray*}\label{Bertrand-curvature}
\overline{\ell}_1(t) &=& 0 ,\\
\overline{\ell}_2(t) &=& (\ell_1(t) \sin \psi(t)+\ell_2(t)\cos \psi(t))\cos \theta(t)-\ell_3(t)\sin \theta(t),\\
\overline{\ell}_3(t) &=& (\ell_1(t) \sin \psi(t)+\ell_2(t)\cos \psi(t))\sin \theta(t)+\ell_3(t)\cos \theta(t),\\
\overline{\ell}_4(t) &=& (\ell_4(t)-\dot{\psi}(t))\cos \theta(t)-(\ell_5(t)\cos \psi(t)-\ell_6(t)\sin \psi(t))\sin \theta(t),\\
\overline{\ell}_5(t) &=& (\ell_4(t)-\dot{\psi}(t))\sin \theta(t)+(\ell_5(t)\cos \psi(t)-\ell_6(t)\sin \psi(t))\cos \theta(t),\\
\overline{\ell}_6(t) &=& -\dot{\theta}(t)+\ell_5(t)\sin \psi(t)+\ell_6(t)\cos \psi(t),\\
\overline{\alpha}(t) &=& \lambda (\ell_1(t)\sin\psi(t)+\ell_2(t)\cos \psi(t))\sin \theta(t)+(\alpha(t)+\lambda \ell_3(t))\cos \theta(t).
\end{eqnarray*}
\end{corollary}
\begin{remark}{\rm
Suppose that $(\gamma,\nu_1,\nu_2,\nu_3)$ and $(\overline{\gamma},\overline{\nu}_1,\overline{\nu}_2,\overline{\nu}_3):I \to \R^4 \times \Delta$ are Bertrand mates which is given by \eqref{Bertrand-5}, \eqref{Bertrand-6}, \eqref{Bertrand-7} and \eqref{Bertrand-8}.
By equation \eqref{Bertrand-condition2}, we have $\overline{\ell}_2(t)=\alpha(t)\sin \theta(t)/\lambda$.
}
\end{remark}

\begin{example}{\rm
Let $(\gamma,\nu_1,\nu_2,\nu_3):[0,2\pi) \to \R^4 \times \Delta$,
\begin{eqnarray*}
\gamma(t)&=&\left(t\sin t+\cos t, -t\cos t+\sin t, t\sin 2t+\frac{1}{2}\cos 2t, -t\cos 2t+\frac{1}{2}\sin 2t\right),\\
\nu_1(t)&=&(-\sin t,\cos t,0,0),\\
\nu_2(t)&=&(0,0,-\sin 2t,\cos 2t),\\
\nu_3(t)&=&\frac{2}{\sqrt{5}}\left(\cos t,\sin t,-\frac{1}{2}\cos 2t,-\frac{1}{2}\sin 2t\right).
\end{eqnarray*}
Note that $t=0$ is a singular point of $\gamma$.
By a direct calculation, $(\gamma,\nu_1,\nu_2,\nu_3)$ is a framed curve.
Then
$$
\bmu(t)=\nu_1(t)\times \nu_2(t)\times \nu_3(t)=\frac{1}{\sqrt{5}}(\cos t,\sin t,2\cos 2t,2\sin 2t)
$$
and the curvature of $(\gamma,\nu_1,\nu_2,\nu_3)$ is given by
$$
(\ell_1(t),\dots,\ell_6(t),\alpha(t))=\left(0, \ -\frac{2}{\sqrt{5}}, \ -\frac{1}{\sqrt{5}}, \ \frac{2}{\sqrt{5}}, \ -\frac{4}{\sqrt{5}}, \ 0,\ \sqrt{5}t\right).$$
If we take $\psi(t)=0$, $\theta(t)$ satisfy the equation $2\cos\theta(t)-(\sqrt{5}t+1)\sin\theta(t)=0$ and $\lambda=-\sqrt{5}$, then
equations \eqref{Bertrand-condition1} and \eqref{Bertrand-condition2} are satisfied.
Hence, $(\gamma,\nu_1,\nu_2,\nu_3)$ is a Bertrand curve by Theorem \ref{framed-Bertrand-equivalent}.
In fact, by equations \eqref{Bertrand-5}, \eqref{Bertrand-6}, \eqref{Bertrand-7} and \eqref{Bertrand-8},
$(\overline{\gamma},\overline{\nu}_1,\overline{\nu}_2,\overline{\nu}_3): [0,2\pi) \to \R^4 \times \Delta$,
\begin{eqnarray*}
\overline{\gamma}(t)&=&\gamma(t)-\sqrt{5}\nu_1(t)\\
&=&\left((t+\sqrt{5})\sin t+\cos t, -(t+\sqrt{5})\cos t+\sin t, t\sin 2t+\frac{1}{2}\cos 2t, -t\cos 2t+\frac{1}{2}\sin 2t\right),\\
\overline{\nu}_1(t)&=&\nu_1(t),\\
\overline{\nu}_2(t)&=&\nu_2(t), \\
\overline{\nu}_3(t)&=&\cos\theta(t)\nu_3(t)-\sin\theta(t)\bmu(t)\\
&=&\frac{1}{\sqrt{(\sqrt{5}t+1)^2+4}}\left(2t\cos t,2t\sin t,-(t+\sqrt{5})\cos 2t,-(t+\sqrt{5})\sin 2t\right)
\end{eqnarray*}
is a framed curve.
Hence, $(\gamma,\nu_1,\nu_2,\nu_3)$ and $(\overline{\gamma},\overline{\nu}_1,\overline{\nu}_2,\overline{\nu}_3)$ are Bertrand metes.
}
\end{example}

\section{Mannheim curves of framed curves}

Let $(\gamma,\nu_1,\nu_2,\nu_3)$ and $(\overline{\gamma},\overline{\nu}_1,\overline{\nu}_2,\overline{\nu}_3):I \to \R^4 \times \Delta$ be framed curves with the curvatures $(\ell_1,\dots,\ell_6,\alpha)$ and $(\overline{\ell}_1,\dots,\overline{\ell}_6,\overline{\alpha})$, respectively.
Suppose that $\gamma$ and $\overline{\gamma}$ are different curves, that is, $\gamma \not\equiv \overline{\gamma}$.

\begin{definition}\label{Mannheim.framed}{\rm
$(1)$ We say that framed curves $(\gamma,\nu_1,\nu_2,\nu_3)$ and $(\overline{\gamma},\overline{\nu}_1,\overline{\nu}_2,\overline{\nu}_3)$ are {\it second type of Mannheim mates} (briefly, {\it second mates} or $(\nu_1,\overline{\nu}_2)$-{\it mates}) if there exists a smooth function $\lambda:I \to \R$ such that $\overline{\gamma}(t)=\gamma(t)+\lambda(t)\nu_1(t)$ and $\nu_1(t)=\overline{\nu}_2(t)$ for all $t \in I$.
\par
$(2)$ We say that framed curves $(\gamma,\nu_1,\nu_2,\nu_3)$ and $(\overline{\gamma},\overline{\nu}_1,\overline{\nu}_2,\overline{\nu}_3)$ are {\it third type of Mannheim mates} (briefly, {\it third mates} or $(\nu_1,\overline{\nu}_3)$-{\it mates}) if there exists a smooth function $\lambda:I \to \R$ such that $\overline{\gamma}(t)=\gamma(t)+\lambda(t)\nu_1(t)$ and $\nu_1(t)=\overline{\nu}_3(t)$ for all $t \in I$.
\par
We also say that $(\gamma,\nu_1,\nu_2,\nu_3):I \to \R^4 \times \Delta$ is a {\it second type of Mannheim curve} (respectively, {\it third type of Mannheim curve}) if there exists a framed curve  $(\overline{\gamma},\overline{\nu}_1,\overline{\nu}_2,\overline{\nu}_3):I \to \R^4 \times \Delta$ such that $(\gamma,\nu_1,\nu_2,\nu_3)$ and $(\overline{\gamma},\overline{\nu}_1,\overline{\nu}_2,\overline{\nu}_3)$ are second mates (respectively, third mates).
}
\end{definition}

\begin{lemma}\label{lambda-const2}
Under the notations in Definition \ref{Mannheim.framed},
if $(\gamma,\nu_1,\nu_2,\nu_3)$ and $(\overline{\gamma},\overline{\nu}_1,\overline{\nu}_2,\overline{\nu}_3)$ are second or third mates, then $\lambda$ is a non-zero constant.
\end{lemma}
\demo
By differentiating $\overline{\gamma}(t)=\gamma(t)+\lambda(t)\nu_1(t)$, we have
$$
\overline{\alpha}(t) \overline{\bmu}(t)
=(\alpha(t)+\lambda(t)\ell_3(t))\bmu(t)+\dot{\lambda}(t)\nu_1(t)+\lambda(t)\ell_1(t)\nu_2(t)+\lambda(t)\ell_2(t)\nu_3(t)
$$
for all $t \in I$.
Since $\nu_1(t)=\overline{\nu}_2(t)$ or $\nu_1(t)=\overline{\nu}_3(t)$, we have $\dot{\lambda}(t)=0$ for all $t \in I$.
Therefore $\lambda$ is a constant.
If $\lambda=0$, then $\overline{\gamma}(t)=\gamma(t)$ for all $t \in I$.
Hence, $\lambda$ is a non-zero constant.
\enD

By definitions of the Bertrand curve, second type and third type of Mannheim curve of a framed curve, we have the following result.
\begin{theorem}\label{relations}
Let $(\gamma,\nu_1,\nu_2,\nu_3):I \to \R^4 \times \Delta$ be a framed curve.
Then the following are equivalent:
\par
$(1)$ $(\gamma,\nu_1,\nu_2,\nu_3)$ is a Bertrand curve.
\par
$(2)$ $(\gamma,\nu_1,\nu_2,\nu_3)$ is a second type of Mannheim curve.
\par
$(3)$ $(\gamma,\nu_1,\nu_2,\nu_3)$ is a third type of Mannheim curve.
\end{theorem}
\demo
Suppose that $(\gamma,\nu_1,\nu_2,\nu_3)$ is a Bertrand curve, that is, there exists a framed curve
$(\overline{\gamma},\overline{\nu}_1,\overline{\nu}_2,\overline{\nu}_3)$ such that $\overline{\gamma}(t)=\gamma(t)+\lambda \nu_1(t)$ and $\nu_1(t)=\overline{\nu}_1(t)$ for all $t \in I$, where $\lambda$ is a non-zero constant.
If we denote $(\widetilde{\nu}_1,\widetilde{\nu}_2,\widetilde{\nu}_3)=(\overline{\nu}_3,\overline{\nu}_1,\overline{\nu}_2)$, then $(\overline{\gamma},\widetilde{\nu}_1,\widetilde{\nu}_2,\widetilde{\nu}_3)$ is also a framed curve with $\nu_1(t)=\widetilde{\nu}_2(t)$ for all $t \in I$.
Hence, $(\gamma,\nu_1,\nu_2,\nu_3)$ is a second type of Mannheim curve and vice versa.
Moreover, if we consider $(\widetilde{\nu}_1,\widetilde{\nu}_2,\widetilde{\nu}_3)=(\overline{\nu}_2,\overline{\nu}_3,\overline{\nu}_1)$, then $(\overline{\gamma},\widetilde{\nu}_1,\widetilde{\nu}_2,\widetilde{\nu}_3)$ is also a framed curve with $\nu_1(t)=\widetilde{\nu}_3(t)$ for all $t \in I$.
Hence, $(\gamma,\nu_1,\nu_2,\nu_3)$ is a third type of Mannheim curve and vice versa.
\enD
\begin{remark}{\rm
We can also prove Theorem \ref{relations} by using direct calculations of existence conditions of second and third types of Mannheim curves.
See \cite{Honda-Takahashi2} for the case of $\R^3$.
}
\end{remark}

\section{Revisit to Bertrand curves of regular space curves in $\R^4$}

Let $\gamma:I \to \R^4$ be a regular space curve with non-degenerate and the arc-length parameter.
Then we have the moving frame $\{\bt(s),\bn_1(s),\bn_2(s),\bn_3(s)\}$ of $\gamma(s)$.
We consider a framed curve by $(\gamma,\bn_1,\bn_2,\bn_3):I \to \R^4 \times \Delta$.
Note that $\bn_1(s) \times \bn_2(s) \times \bn_3(s)=-\bt(s)$.
By the Frenet-Serret formula (see \S 2), we have
$$
\left(
\begin{array}{c}
\bn'_1(s)\\
\bn'_2(s)\\
\bn'_3(s)\\
-\bt'(s)\\
\end{array}
\right)
=
\left(
\begin{array}{cccc}
0&\kappa_2(s)&0&\kappa_1(s)\\
-\kappa_2(s)&0&\kappa_3(s)&0\\
0&-\kappa_3(s)&0&0\\
-\kappa_1(s)&0&0&0
\end{array}
\right)
\left(
\begin{array}{c}
\bn_1(s)\\
\bn_2(s)\\
\bn_3(s)\\
-\bt(s)
\end{array}
\right), \gamma'(s)=-(-\bt(s)).
$$
Thus,
\begin{eqnarray*}
\ell_1(s)=\kappa_2(s), \ \ell_2(s)=0, \ \ell_3(s)=\kappa_1(s), \ \ell_4(s)=\kappa_3(s), \ \ell_5(s)=0, \  \ell_6(s)=0, \ \alpha(s)=-1.
\end{eqnarray*}
By Theorem \ref{framed-Bertrand-equivalent}, equations \eqref{Bertrand-condition1} and \eqref{Bertrand-condition2} are given by
\begin{eqnarray*}
&&\ell_1(s)\cos \psi(s)-\ell_2(s)\sin \psi(s)=\kappa_2(s)\cos \psi(s)=0,\\
&&\lambda(\ell_1(s)\sin \psi(s)+\ell_2(s)\cos \psi(s))\cos \theta(s)-(\alpha(s)+\lambda \ell_3(s))\sin \theta(s)\\
&&=\lambda \kappa_2(s) \sin \psi(s) \cos \theta(s)-(-1+\lambda \kappa_1(s)) \sin \theta(s)=0.
\end{eqnarray*}
If $\psi(s)=\pi/2$ and $\lambda \kappa_2(s)\cos \theta(s)+(1-\lambda \kappa_1(s))\sin \theta(s)=0$ for all $s \in I$, then $(\gamma,\bn_1,\bn_2,\bn_3)$ is a Bertrand curve as a framed curve by Theorem \ref{framed-Bertrand-equivalent}.
\par
Next, we consider Bertrand curves of regular space curves which have more geometric meaning as compared to Bertrand curves of framed curves.
It seems that the non-degenerate condition is a strong assumption for getting a moving frame.
We consider mild assumptions between regular space curves with non-degenerate condition and framed curves.
This phenomena occur higher dimensional cases.
It cannot consider regular space curves in 3-dimensional space.
\par
Let $\gamma:I \to \R^4$ be a regular space curve.
In order to define Bertrand curves in the sense of regular space curves in \S 2, we only need the first normal vector $\bn_1$ of $\gamma$.
Therefore, we assume that
$$
|\dot{\gamma}(t)|^2|\ddot{\gamma}(t)|^2-(\dot{\gamma}(t) \cdot \ddot{\gamma}(t))^2 >0
$$
for all $t \in I$.
Namely, $\kappa_1(t) >0$ for all $t \in I$.
Moreover, suppose that there exists a smooth mapping $\nu_2: I \to S^3$ such that $(\gamma,\bn_1,\nu_2,\nu_3):I \to \R^4 \times \Delta$ is a framed curve, where $\bt \times \bn_1 \times \nu_2=\nu_3$.

Let $\gamma$ and $\overline{\gamma}:I \to \R^4$ be regular space curves with the first curvatures $\kappa_1$ and $\overline{\kappa}_1$ are positive.
Suppose that there exist smooth mappings $\nu_2$ and $\overline{\nu}_2:I \to S^3$ such that $(\gamma,\bn_1,\nu_2,\nu_3)$ and $(\overline{\gamma},\overline{\bn}_1,\overline{\nu}_2,\overline{\nu}_3):I \to \R^4 \times \Delta$ are framed curves, where $\bt \times \bn_1 \times \nu_2=\nu_3$ and $\overline{\bt} \times \overline{\bn}_1 \times \overline{\nu}_2=\overline{\nu}_3$.
\begin{definition}\label{Bertrand-regular-framed}{\rm
We say that the above framed curves $(\gamma,\bn_1,\nu_2,\nu_3)$ and $(\overline{\gamma},\overline{\bn}_1,\overline{\nu}_2,\overline{\nu}_3):I \to \R^4 \times \Delta$ are {\it Bertrand mates} if there exists a smooth function $\lambda:I \to \R$ such that $\overline{\gamma}(t)=\gamma(t)+\lambda(t)\bn_1(t)$ and $\bn_1(t)=\pm \overline{\bn}_1(t)$ for all $t \in I$.
We also say that $(\gamma,\bn_1,\nu_2,\nu_3):I \to \R^4 \times \Delta$ is a {\it Bertrand curve} if there exists a framed curve  $(\overline{\gamma},\overline{\bn}_1,\overline{\nu}_2,\overline{\nu}_3):I \to \R^4 \times \Delta$ such that $(\gamma,\bn_1,\nu_2,\nu_3)$ and $(\overline{\gamma},\overline{\bn}_1,\overline{\nu}_2,\overline{\nu}_3)$ are Bertrand mates.
}
\end{definition}
Since $(\gamma,\bn_1,\nu_2,\nu_3)$ is a special case of a framed curve, we have the following results, see in \S \ref{Bertrand-framed}.
\begin{corollary}
Under the notations in Definition \ref{Bertrand-regular-framed},
if $(\gamma,\bn_1,\nu_2,\nu_3)$ and $(\overline{\gamma},\overline{\bn}_1,\overline{\nu}_2,\overline{\nu}_3)$ are Bertrand mates, then $\lambda$ is a non-zero constant.
\end{corollary}
Suppose that $s$ is the arc-length parameter of $\gamma$.
For a framed curve $(\gamma,\bn_1,\nu_2,\nu_3): I \to \R^4 \times \Delta$, we have $\bn_1 \times \nu_2 \times \nu_3=-\bt$ and
 $$
\left(
\begin{array}{c}
\bn'_1(s)\\
\nu'_2(s)\\
\nu'_3(s)\\
-\bt'(s)\\
\end{array}
\right)
=
\left(
\begin{array}{cccc}
0&\ell_1(s)&\ell_2(s)&\kappa_1(s)\\
-\ell_1(s)&0&\ell_3(s)&0\\
-\ell_2(s)&-\ell_3(s)&0&0\\
-\kappa_1(s)&0&0&0
\end{array}
\right)
\left(
\begin{array}{c}
\bn_1(s)\\
\nu_2(s)\\
\nu_3(s)\\
-\bt(s)
\end{array}
\right), \gamma'(s)=-(-\bt(s)).
$$
Thus, the curvature of the framed curve $(\gamma,\bn_1,\nu_2,\nu_3)$ is given by $(\ell_1,\ell_2,\kappa_1,\ell_3,0,0,-1)$.

We give a necessary and sufficient condition of a Bertrand curve in the sense of Definition \ref{Bertrand-regular-framed}.

\begin{theorem}\label{regular-framed-Bertrand-equivalent}
Under the notations in Definition \ref{Bertrand-regular-framed}, suppose that $(\gamma,\bn_1,\nu_2,\nu_3)$ is a framed curve with the curvature $(\ell_1,\ell_2,\kappa_1,\ell_3,0,0,-1)$, $s$ is the arc-length parameter of $\gamma$ and $\lambda$ is a non-zero constant.
Then $(\gamma,\bn_1,\nu_2,\nu_3)$ and $(\overline{\gamma},\overline{\bn}_1,\overline{\nu}_2,\overline{\nu}_3)$ are Bertrand mates with $\overline{\gamma}(s)=\gamma(s)+\lambda \bn_1(s)$ if and only if there exist smooth functions $\psi,\theta:I \to \R$ such that
\begin{eqnarray}
&&\ell_1(s)\cos \psi(s)-\ell_2(s)\sin \psi(s)=0, \label{regular-Bertrand-condition1}\\
&&\lambda(\ell_1(s) \sin \psi(s)+ \ell_2(s)\cos \psi(s))\cos \theta(s)+(1-\lambda \kappa_1(s))\sin \theta(s)=0, \label{regular-Bertrand-condition2}\\
&&h(s) = (1-\lambda\kappa_1(s))^2((1-\lambda \kappa_1(s))\kappa_1(s)-\ell_1^2(s)-\ell_2^2(s))^2 \nonumber\\
&&\quad +(1-\lambda\kappa_1(s))(\ell'_1(s)-\ell_2(s)\ell_3(s))\left((1-\lambda \kappa_1(s))(\ell_1'(s)-\ell_2(s)\ell_3(s))+2\lambda\kappa_1'(s)\ell_1(s)\right) \nonumber\\
&&\quad +(1-\lambda\kappa_1(s))(\ell'_2(s)-\ell_1(s)\ell_3(s))\left((1-\lambda \kappa_1(s))(\ell_2'(s)+\ell_1(s)\ell_3(s))+2\lambda\kappa_1'(s)\ell_2(s)\right) \nonumber\\
&&\quad +(\lambda^2\kappa_1'(s)^2+((1-\lambda\kappa_1(s))\kappa_1(s)-\ell_1^2(s)-\ell_2^2(s))^2)(\ell_1^2(s)+\ell_2^2(s)) \nonumber\\
&&\quad +(\ell_2(s)(\ell_1'(s)-\ell_2(s)\ell_3(s))-\ell_1(s)(\ell_2'(s)+\ell_1(s)\ell_3(s)))^2 >0 \label{regular-Bertrand-condition3}
\end{eqnarray}
for all $s \in I$.
\end{theorem}
\demo
By Theorem \ref{framed-Bertrand-equivalent},
replace the curvature $(\ell_1,\dots,\ell_6,\alpha)$ to $(\ell_1,\ell_2,\kappa_1,\ell_3,0,0,-1)$ in equations \eqref{Bertrand-condition1} and \eqref{Bertrand-condition2}, we have conditions \eqref{regular-Bertrand-condition1} and \eqref{regular-Bertrand-condition2}.
Hence, it is enough to give the condition that the first curvature $\overline{\kappa}_1$ of $\overline{\gamma}$ is positive.
Since
\begin{eqnarray*}
\dot{\overline{\gamma}}(s) &=& (1-\lambda \kappa_1(s))\bt(s)+\ell_1(s)\nu_2(s)+\ell_2(s)\nu_3(s),\\
\ddot{\overline{\gamma}}(s) &=& -\lambda {\kappa}_1'(s)\bt(s)+((1-\lambda\kappa_1(s))\kappa_1(s)-\ell_1^2(s)-\ell_2^2(s))\bn_1(s)\\
&&+(\ell_1'(s)-\ell_2(s)\ell_3(s))\nu_2(s)+(\ell_2'(s)+\ell_1(s)\ell_3(s))\nu_3(s),
\end{eqnarray*}
we have
\begin{eqnarray*}
|\dot{\overline{\gamma}}(s)|^2 &=& (1-\lambda \kappa_1(s))^2+\ell_1^2(s)+\ell^2_2(s),\\
|\ddot{\overline{\gamma}}(s)|^2 &=& \lambda^2 \kappa_1'(s)^2+((1-\lambda\kappa_1(s))\kappa_1(s)-\ell_1^2(s)-\ell_2^2(s))^2\\
&&+(\ell_1'(s)-\ell_2(s)\ell_3(s))^2+(\ell_2'(s)+\ell_1(s)\ell_3(s))^2,\\
\dot{\overline{\gamma}}(s) \cdot \ddot{\overline{\gamma}}(s)&=&-(1-\lambda\kappa_1(s))\lambda \kappa_1'(s)+\ell_1(s)(\ell_1'(s)-\ell_2(s)\ell_3(s))+\ell_2(s)(\ell_2'(s)+\ell_1(s)\ell_3(s)).
\end{eqnarray*}
By a direct calculation,
$|\dot{\overline{\gamma}}(s)|^2|\ddot{\overline{\gamma}}(s)|^2-(\dot{\overline{\gamma}}(s) \cdot \ddot{\overline{\gamma}}(s))^2>0$ is equivalent to the condition $h(s)>0$ for all $s \in I$.
\enD
\begin{remark}{\rm
Under the same assumptions in Theorem \ref{regular-framed-Bertrand-equivalent}, $(\gamma,\bn_1,\nu_2,\nu_3)$ and $(\overline{\gamma},\overline{\bn}_1,\overline{\nu}_2,\overline{\nu}_3)$ are Bertrand mates with $\overline{\gamma}(s)=\gamma(s)+\lambda \bn_1(s)$.
Then the first curvature $\overline{\kappa}_1(s)$ of $\overline{\gamma}(s)$ is given by
$\sqrt{h(s)}/\sqrt{(1-\lambda \kappa_1(s))^2+\ell_1^2(s)+\ell_2^2(s)}$.
}
\end{remark}

\begin{example}{\rm
Let $\gamma:[0,2\pi) \to \R^4$ be a regular space curve which is given by
$$
\gamma(s)=\frac{1}{\sqrt{2}}\left(\frac{1}{a} \cos as, \frac{1}{a} \sin as, \frac{1}{b} \cos bs, \frac{1}{b}\sin bs \right),
$$
where $a$ and $b$ are non-zero constants with $a^2 \not=b^2$.
Since
$$
\gamma'(s)=\frac{1}{\sqrt{2}}(-\sin as, \cos as, -\sin bs, \cos bs),
$$
we have $|\gamma'(s)|=1$ for all $s \in [0,2\pi)$.
Hence $s$ is the arc-length parameter of $\gamma$, $\bt(s)=\gamma'(s)$ and $|\bt'(s)|=\sqrt{a^2+b^2}/\sqrt{2}$.
It follows that $\kappa_1(s)=\sqrt{a^2+b^2}/\sqrt{2}$ and
$$
\bn_1(s)=\frac{\bt'(s)}{|\bt'(s)|}=\frac{1}{\sqrt{a^2+b^2}}(-a \cos as,-a \sin as,-b \cos bs, -b \sin bs).
$$
If we take $\nu_2: [0,2\pi) \to S^3$,
$\nu_2(s)=(1/{\sqrt{2}})(-\sin as, \cos as, \sin bs, -\cos bs)$, then
\begin{eqnarray*}
\nu_3(s) &=& \bt(s) \times \bn_1(s) \times \nu_2(s) \\
&=& \frac{1}{\sqrt{a^2+b^2}}(b\cos as, b \sin as, -a \cos bs, -a \sin bs).
\end{eqnarray*}
It follows that $(\gamma,\bn_1,\nu_2,\nu_3): [0,2\pi) \to \R^4 \times \Delta$ is a framed curve.
By a direct calculation, the curvature of $(\gamma,\bn_1,\nu_2,\nu_3)$ is given by
$$
(\ell_1(t),\ell_2(t),\kappa_1(t),\ell_3(t),0,0,-1)=\left(\frac{-a^2+b^2}{\sqrt{2(a^2+b^2)}},0,\frac{\sqrt{a^2+b^2}}{\sqrt{2}},-\frac{\sqrt{2}ab}{\sqrt{a^2+b^2}},0,0,-1\right).
$$
If we take $\psi(s)=\pi/2$, then equation \eqref{regular-Bertrand-condition1} is satisfied.
Suppose that $1-\lambda \kappa_1 \not=0$.
Then there exists a constant $\theta$ such that equation \eqref{regular-Bertrand-condition2} is satisfied.
By a direct calculation, we have
$$
h=((1-\lambda \kappa_1)\kappa_1-\ell_1^2)^2(\ell_1^2+(1-\lambda \kappa_1)^2)+(\ell_1\ell_3)^2(\ell_1^2-(1-\lambda \kappa_1)^2).
$$
If $h>0$, then $(\gamma,\bn_1,\nu_2,\nu_3)$ is a Bertrand curve by Theorem \ref{regular-framed-Bertrand-equivalent}.
\par
As a concrete example, we consider $a=1, b=\sqrt{3}$ and $\lambda=\sqrt{2}/4$.
In this case, $\psi(s)=\pi/2, 1-\lambda \kappa_1=1/2 \not=0, \theta={\rm Tan^{-1}}(-1/2)$ and $h=3(2-\sqrt{2})/8>0$.
It follows that $(\gamma,\bn_1,\nu_2,\nu_3)$ is a Bertrand curve.
}
\end{example}

Shun'ichi Honda,
\par
Chitose Institute of Science and Technology, Chitose 066-8655, Japan,
\par
E-mail address: s-honda@photon.chitose.ac.jp
\\

Masatomo Takahashi,
\par
Muroran Institute of Technology, Muroran 050-8585, Japan,
\par
E-mail address: masatomo@mmm.muroran-it.ac.jp
\\

Haiou Yu,
\par
Jilin University of Finance and Economic, Changchun 130117, China,
\par
E-mail address: yuhaiou@jlufe.edu.cn

\end{document}